\documentclass[11pt]{amsart}
\usepackage{amsmath,amssymb,bm}
\usepackage{verbatim,a4wide} 
\newtheorem{theorem}{Theorem}[section]
\newtheorem{lemma}[theorem]{Lemma}
\newtheorem{proposition}[theorem]{Proposition}

\newtheorem{claim}{Claim}[section]
\theoremstyle{remark}

\numberwithin{equation}{section}
\newcommand{\R}{\mathbb{R}}

\newcommand{\N}{\mathbb{N}}
\newcommand{\Z}{\mathbb{Z}}
\newcommand{\T}{\mathbb{T}}

\newcommand{\pd}{\partial}
\newcommand{\eps}{\varepsilon}
\DeclareMathOperator{\sech}{sech}

\DeclareMathOperator{\spann}{span}

\begin{document}
\title
[Stability of the line soliton of the KP-II equation]
{Stability of the line soliton of the KP-II equation under periodic transverse perturbations}
\author{Tetsu Mizumachi and Nikolay Tzvetkov}
\address{Department of Mathematics, Kyushu University, Fukuoka, 819-0395, Japan and University of 
Cergy-Pontoise, UMR CNRS 8088, Cergy-Pontoise, F-95000}
\begin{abstract}
We prove the nonlinear stability of the KdV solitary waves considered as solutions of
the KP-II equation, with respect to periodic transverse perturbations.
\end{abstract}
\maketitle
\section{Introduction}
\label{sec:intro}
Our goal here is to prove the nonlinear stability of the KdV solitary waves
considered as solutions of the Kadomtsev-Petviashvili-II (KP-II) equation
\begin{equation}\label{eq:KPII}
\partial_x(\pd_tu+\pd_x^3u+3\pd_x(u^2))+3\partial_y^2u=0
\end{equation}
with respect to periodic transverse perturbations.
In this paper, we consider \eqref{eq:KPII} for $(x,y)\in \R\times\T$, 
where $\T=\R/(2\pi\Z)$ denotes the one dimensional torus.
 
The well-posedness of \eqref{eq:KPII} is studied in \cite{MST}, where it is shown that \eqref{eq:KPII} is globally well-posed for initial data in
$H^s(\R_x\times\T_y)$, $s\geq 0$.
Roughly speaking, it is shown that for every $u_0\in H^s(\R\times \T)$ with
$s\geq 0$, there is a unique solution of \eqref{eq:KPII}
which belongs to $C(\R;H^s(\R\times \T))$.
Moreover the flow map is continuous (and even uniformly continuous on
bounded sets) in the phase space $H^s(\R\times \T)$.
The proof is based on the ideas introduced in the purely periodic case 
in the work of Bourgain~\cite{Bourgain}.
For other contributions on the Cauchy problem of the KP-II equation with
different spatial domains, we refer to \cite{GPS,Hadac,HHK,IM,Tak,TT,Tz}.

Let us now turn to the stability questions.
Let
$$\varphi_c(x)\equiv c\cosh^{-2}\Big(\sqrt{\frac{c}{2}}\,x\Big),\quad c>0.
$$
Then $\varphi_c(x-2ct)$ is a solitary wave solution of the KdV equation and
also a solution of \eqref{eq:KPII}. It is well-known that $\varphi_c(x-2ct)$ is
orbitally stable as a solution of the KdV equation (see \cite{Benjamin,Bona}).
Our goal here is to show that $\varphi_c(x-2ct)$ remains stable as a solution
of the KP-II equation subject to perturbations which are periodic in the
transversal direction. Now let us introduce our result.
\begin{theorem}\label{thm1}
For every $\varepsilon>0$, there exists a $\delta>0$ such that
if the initial data of \eqref{eq:KPII} satisfies
$\|u_0-\varphi_c\|_{L^2(\R_x\times\T_y)}< \delta$,
the corresponding solution of \eqref{eq:KPII} satisfies
\begin{equation*}
\inf_{\gamma\in\R}\|u(t,x,y)-\varphi_c(x+\gamma)\|_{L^2(\R_x\times\T_y)}<\varepsilon,\quad \forall\, t\in\R.
\end{equation*}
Moreover, there exists a constant $\tilde{c}$ satisfying $\tilde{c}-c=O(\delta)$ and a modulation parameter $x(t)$ satisfying $\limsup_{t\to\infty}|x(t)/t-2c|=O(\delta)$ and such that
\begin{equation}\label{eq:asympst}
\lim_{t\to\infty}\|u(t,x,y)-\varphi_{\tilde{c}}(x-x(t))\|_{L^2((x\ge ct)\times\T_y)}=0.
\end{equation}
\end{theorem}

This result confirms the heuristic analysis of the Kadomtsev-Petviashvili
seminal paper \cite{KP}. 
In the recent work \cite{VA}, a stability result for the KP-II line soliton is studied
by the inverse scattering method.
There are several differences between \cite{VA} and Theorem~\ref{thm1}.
For instance in \cite{VA}, localized perturbations belonging to
weighted spaces are considered and thus one does not need to involve
the modulation parameter $x(t)$ in the stability statement.  On the
other hand, for the periodic perturbations considered in this paper,
the modulation by translation is unavoidable (as it is for the KdV
equation).  In addition, our approach is apparently less dependent on
the integrability features of the KP-II equation.
We refer to \cite{Cu} for a recent study on transverse stability for
Hamiltonian PDE's. 

Let us now explain the main ideas and difficulties behind the proof of Theorem~\ref{thm1}.
The problem \eqref{eq:KPII} has a Lax pair structure (see \cite{ZS}) and thus it has, at least formally,
an infinite sequence of conservation laws. Unfortunately, these conservation laws do not seem easy to use
for dynamical issues. The only conservation law for the KP-II equation used in this paper is the $L^2$
conservation law, i.e. the quantity 
$$
N(u(t,\cdot))=\int_{\R_x\times\T_y} u^2(t,x,y)dxdy
$$ 
is conserved (independent of $t$) by the flow of \eqref{eq:KPII} established
by the well-posedness result in \cite{MST}.
Indeed, the first two terms of the Hamiltonian
$$
E(u(t,\cdot))=
\int_{\R_x\times\T_y}
\big(u^2_x(t,x,y)-3(\pd_x^{-1}\pd_yu(t,x,y))^2-2u^3(t,x,y)\big)dxdy,
$$
which is one of the conserved quantities of \eqref{eq:KPII},
have the opposite sign. The infinite dimensional indefiniteness is a
serious obstruction to use the Hamiltonian in controlling the long
time behavior of the KP-II equation. 
In particular, we cannot use the standard approach to prove stability based on
the fact that the line soliton $\varphi_c(x-2ct)$ is a minimizer of
the functional $E(u)$ on the manifold $\{u\in H^1 \,|\, N(u)=N(\varphi_c)\}$.

In this paper, we aim to prove modulational stability in the $L^2$-framework.
For that purpose, we follow the idea of the work by Merle and Vega \cite{MV}
which prove orbital/asymptotic stability of KdV 1-soliton in $L^2$.
The idea of Merle and Vega \cite{MV} is to lift up a solution around
a 1-soliton of KdV to a solution
around a kink solution of the modified KdV equation by using
the Miura transform. Since a kink is not in the energy class, the Miura
transform eliminates the scaling  freedom which generates the only direction we must be afraid of,
to argue stability by using the $L^2$-conservation law. 
The other merit of using the Miura transform is that it gains 1 more derivative and makes it possible to argue
stability of kinks by a standard energy method.
The Miura transform is one of the B\"acklund transformations that enable
us to observe behavior of solutions in a \lq\lq simpler coordinate\rq\rq.
For example, it enables us to prove linear stability of solitons in
a simple way (see \cite{MP}). 

The Miura transform associated to the KP-II equation \eqref{eq:KPII} is a heat
operator (see \cite{W,KMa}).
The main point is that it has a  similar structure near the solitary wave 
with the Miura transform of the KdV equation which makes Merle and Vega's
approach applicable.
Here we need to replace the ODE argument of \cite{MV} by a suitable Lyapunov functional argument.
Once the crucial analysis of the Miura transform near a solitary wave is performed, we can
argue stability of line solitons of \eqref{eq:KPII} through
stability of kink solutions of the mKP-II equation (the equation obtained from \eqref{eq:KPII}
after applying the Miura transform). This essentially explains our approach.

To obtain the asymptotic stability result, we use monotonicity coming from a 
Kato type smoothing effect (see \cite{dBM,KMa}). 
Our proof is  simpler than a paper by Martel and Merle\ \cite{MM}
because thanks to the Miura transform, we only need monotonicity property of
small solutions to the KP-II equation.
In fact, we do not need the modulation equation for the amplitude of the
main line soliton because it is {\it \`a priori} determined through
the Miura transform.

It would be interesting to extend our results to fully localized perturbations (belonging to $H^s(\R^2)$) of the  KdV soliton under the KP-II flow. In this case the study of the linearization of the Miura transform $M^c_{+}$ in a neighborhood of $Q_c$ is much more delicate since the ODE analysis degenerate when the transverse frequencies tends to zero. On the other hand if one succeeds to resolve this issue then it would become possible to use the critical space analysis of \cite{HHK} to get a linear behavior of the solution at the left of the solitary wave and thus give a more precise asymptotic stability statement. We plan to study this phenomenon elsewhere.

The situation changes radically if we replace  \eqref{eq:KPII} by the KP-I equation.
\begin{equation}\label{kp1}
\partial_x(\pd_tu+\pd_x^3u+3\pd_x(u^2))-3\partial_y^2u=0.
\end{equation}
Indeed, it is known since the work of Zakharov \cite{Z} that $\varphi_c(x-2ct)$ is unstable as a solution of \eqref{kp1}.
The proof of Zakharov is based on the integrability features of \eqref{kp1}. We refer to the recent works \cite{RT1,RT2}
for proofs of the instability of $\varphi_c(x-2ct)$ as a solution of \eqref{kp1} independent of the integrability.
These proofs have the advantage to apply to more involved Hamiltonian models.
Let us also refer to \cite{IKT} and the references therein for the quite intricate issues around the well-posedness in Sobolev spaces of \eqref{kp1}.  

Let us complete this introduction by a remark concerning possible extensions.
For some bidirectional model equations such as the FPU lattice equation,
the Hamiltonian is the only useful conservation law as $L^2$-norm is for the
KP-II equation.  Friesecke and Pego \cite{FP} and Mizumachi \cite{Mi1}
prove stability of solitary waves using strong linear stability of solitary
waves in a weighted space.
However, for PDEs such as the water wave models, their approach could require
smallness of higher order 
Sobolev norms that does not follow from conservation laws.
The method of Merle and Vega we use in this paper would suggest that a
lifting argument could help to handle such difficulties in a stability analysis.
We believe that this is an interesting line for further research.
\section{The mKP-II equation}
\label{sec:onmkp2}
This section is devoted to the well-posedness of the mKP-II equation posed on $\R\times \T$.
Here we will strongly rely on the arguments of a recent work by Kenig and 
Martel \cite{KMa}.

We start by defining the anti-derivative operator $\pd_x^{-1}$ via the Fourier transform for
functions $u\in L^2(\R\times \T)$ such that
$\xi^{-1}\hat{u}(\xi,n)\in L^2(\R\times \Z)$ (where $\hat{u}$ denotes the
Fourier transform of $u$). Namely
$
\pd_x^{-1}u=\mathcal{F}_{\xi,n}^{-1}
((i\xi)^{-1}\hat{u}(\xi,n)),
$
where $\mathcal{F}_{\xi,n}^{-1}$ is the inverse Fourier transform.

Let us remark that one may also consider the ``integrated'' form of \eqref{eq:KPII}, namely
\begin{equation}\label{KPII_integrated}
\pd_tu+\pd_x^3u+3\pd_x^{-1}\pd_y^2u+3\pd_x(u^2)=0\,,\quad
u(0,x,y)=u_{0}(x,y).
\end{equation}
The equation \eqref{KPII_integrated} is of the first order in $t$ but one needs
to define $\pd_x^{-1}\pd_y^2u$.
Since the nonlinearity is differentiated with respect to $x$,
this problem only concerns the free evolution.
Thanks to Lemma~1 in \cite{KMa}, one may define $\pd_x^{-1}\pd_y$ of
the free evolution as an $L^1_{loc}(\R^3)$ function even for data
which does not satisfy a constraint $\int u=0$, for example only in $L^2$
(we refer to \cite{MST_contr} for further results in this direction).

We next introduce the space $\mathcal{E}(\R_x\times\T_y)$ which will play an important role 
in the analysis. Let 
$$
\mathcal{E}(\R_x\times\T_y)
=\{u\in L^2(\R\times \T)\,\,:\,\,\|u\|_{\mathcal{E}(\R_x\times\T_y)}<\infty\},
$$
where
\begin{eqnarray*}
\|u\|_{\mathcal{E}(\R_x\times\T_y)}^2
& = & \sum_{n\in\Z}\int_\R
(1+\xi^2+\xi^{-2}n^2)|\hat{u}(\xi,n)|^2\,d\xi
\\
& = &
\|u\|_{L^2}^2+\|\partial_x u\|_{L^2}^2+\|\partial_{x}^{-1}\partial_y u\|_{L^2}^2\,.
\end{eqnarray*}
For each $n\neq 0$, we see that $u_n(x)\equiv\mathcal{F}_\xi^{-1}\hat{u}(\xi,n)$
admits an anti-derivative if $u\in \mathcal{E}$ and that
$(\pd_x^{-1}\pd_y u,v)=(u,\pd_x^{-1}\pd_y v)$ if $u$, $v\in \mathcal{E}$.
Here $(\cdot,\cdot)$ denotes the scalar product of $L^2(\R_x\times\T_y)$ .

We have the following non-isotropic Sobolev inequality for functions in $\mathcal{E}(\R_x\times\T_y)$.
\begin{lemma}\label{Besov}
There exists $C>0$ such that for every $u\in \mathcal{E}(\R_x\times\T_y)$
and $p\in [2,6]$,
$$
\|u\|_{L^p}\leq C\|u\|_{L^2}^{\frac{6-p}{2p}}\|\partial_x u\|_{L^2}^{\frac{p-2}{p}}
\|\partial_x^{-1}\partial_y u\|_{L^2}^{\frac{p-2}{2p}}\,.
$$
\end{lemma}
For a proof of this lemma, we refer to \cite{BIN,Tom} or \cite{MST07} (Lemma~2, page~783). 
The proof on \cite{MST07} is performed
for functions on $\R^2$ but the proof works equally well in the $\R_x\times\T_y$ setting.
In order to motivate the mKP-II equation, we now introduce the Miura transforms
that we use in this paper.  For $c>0$ and $v\in\mathcal{E}$, we set
$$
M^c_{\pm}(v)=\pm \partial_x v+\partial_x^{-1}\partial_{y}v-v^2+\frac{c}{2}\,.
$$
Observe that $M^c_{\pm}$ is invariant by translation, namely
$$
M^c_{\pm}(v(x+\alpha))=(M^c_{\pm}(v))(x+\alpha),\quad\forall\,\alpha\in\R.
$$
Using Lemma~\ref{Besov}, we obtain that if a sequence $\{u_n\}$
converges to a limit $u$ in $\mathcal{E}(\R_x\times\T_y)$ the sequence
$\{M^c_{\pm}(u_n)- M^c_{\pm}(u)\}$ converges to $0$ in $L^2(\R_x\times\T_y)$.

The transformations $M^c_{\pm}$ relate the KP-II equation to the mKP-II
equation (mKP-II)
which reads
\begin{equation}\label{eq:MKPII}
\pd_tv+\pd_x^3v+3\pd_x^{-1}\pd_y^2v-6v^2\pd_xv+6\pd_xv\pd_x^{-1}\pd_yv=0.
\end{equation}
At least formally, if $v(t,x,y)$ is a solution of the mKP-II equation
\eqref{eq:MKPII}, then for $c>0$, $u_{\pm}$ defined by
$$ 
u_{\pm}(t,x,y)\equiv M^c_{\pm}(v)(t,x-3ct,y)
$$
are solutions of the KP-II equation \eqref{eq:KPII}. The last statement can be directly verified
(see e.g. \cite{KMa}, Appendix~A) for
sufficiently smooth solutions in ${\mathcal E}$ and we will only use it in such a situation in this paper.

The line soliton of the KP-II equation is related to the kink $Q_c$ defined by
$$
Q_c(x)=\sqrt{\frac{c}{2}}\tanh\Big(\sqrt{\frac{c}{2}}\,x\Big).
$$
We see that $Q_c(x+ct)$ is a solution of \eqref{eq:MKPII} and moreover
\begin{equation}\label{algebra}
M_{+}^c(Q_c)=\varphi_c,\quad M_{-}^c(Q_c)=0.
\end{equation}
Let $Z=\{u\in H^{8}(\R_x\times\T_y)\,:\, \pd_x^{-1}\pd_yu, \, \pd_xu \in H^8(\R_x\times\T_y)\}.$
In this section, we will prove a global well-posedness result for
\eqref{eq:MKPII} with data \begin{equation}\label{data}
v(0,x,y)=Q_{c}(x)+w_{0}(x,y),\quad w_0\in Z.
\end{equation}
It turns out that one can apply arguments similar to the work by 
Kenig and Martel \cite{KMa} to have the following result. 
\begin{proposition}\label{prop:wpmkp}
For every $w_0\in Z$, there exists a unique global in time solution of
\eqref{eq:MKPII} with data \eqref{data} such that
$$
v(t,x,y)=Q_c(x+ct)+w(t,x,y),\quad w\in C(\R;Z)\,.
$$ 
Moreover $v$ satisfies the conservation laws
\begin{equation}\label{cons_law}
\|M_{\pm}^c(v)(t,x,y)\|_{L^2(\R_x\times\T_y)}=
\|M_{\pm}^c(v)(0,x,y)\|_{L^2(\R_x\times\T_y)}<\infty.
\end{equation}
\begin{proof}
We need to solve the equation
\begin{equation}\label{MKPII_mod}
\pd_t w+\pd_x^3w +3\pd_x^{-1}\pd_y^2 w
-2\partial_x((w+\tilde{Q}_c)^3-\tilde{Q}_c^3)+6\pd_x w\pd_x^{-1}\pd_yw
+6\tilde{Q}_c'\pd_x^{-1}\pd_yw=0 
\end{equation}
with data
\begin{equation}\label{data_mod}
w(0,x,y)=w_{0}(x,y),\quad w_0\in Z,
\end{equation}
where $\tilde{Q}_c\equiv Q_c(x+ct)$.
The construction of local solutions for a regularized version of
\eqref{MKPII_mod}
\begin{equation}\label{eq:regMKPII_mod}
  \begin{split}
&  \pd_t w+\eps\pd_x^4w+\eps^5\pd_y^4w+\pd_x^3w +3\pd_x^{-1}\pd_y^2 w
\\ &
-2\partial_x((w+\tilde{Q}_c)^3-\tilde{Q}_c^3)+6\pd_x w\pd_x^{-1}\pd_yw
+6\tilde{Q}_c'\pd_x^{-1}\pd_yw=0     
  \end{split}
\end{equation}
 can be done as in \cite{KMa}, where 
$v\in\{u\in H^8(\R^2)\,:\, \pd_xu,\, \pd_x^{-1}\pd_yu\in H^8(\R^2)\}$ 
and $\tilde{Q}_c$ is replaced by $0$.  The main point is
a variant of the Kato smoothing effect (Lemma~1 in \cite{KMa}) which
works equally well in the case $\R_x\times\T_y$. Indeed the crucial effect of
the change of variables used in the proof of Lemma~1 in \cite{KMa} is
only in the Fourier variable corresponding to $x$. The arguments used
in the Fourier variable corresponding to $y$ rely only on the Plancherel
theorem and thus the proof of \cite[Lemma~1]{KMa} is transported
directly to the $\R_x\times\T_y$ framework. All other arguments in the
analysis of \eqref{eq:regMKPII_mod} are
independent of the geometry of the spatial domain.  Finally, the
additional terms coming from the presence of $\tilde{Q}_c$ can be
handled similarly to \cite{KMa}.  Indeed $\tilde{Q}_c'\in Z$ which
makes that the term $6\tilde{Q}_c'\pd_x^{-1}\pd_yw$ can be treated
exactly as $6\pd_x w\pd_x^{-1}\pd_yw$.  The new terms coming from the
contributions of $(w+\tilde{Q}_c)^3-\tilde{Q}_c^3$ are
$3\tilde{Q}_c^2w$ and $3\tilde{Q}_cw^2$ and one may readily check that
they do not affect the analysis of \cite[pp.~2462-2463]{KMa}
(they are even slightly easier to handle than $w^3$).

Further one gets local solutions of \eqref{eq:regMKPII_mod}
on a time interval independent of the regularization parameter $\eps$
by a classical compactness argument and a reasoning similar to
the proof of Lemma~\ref{transfert} below,
the argument being even simpler since in Lemma~11 in \cite{KMa},
one needs to establish an energy inequality for \eqref{eq:regMKPII_mod} of the
type $\dot{z}(t)\leq C(1+(z(t))^\gamma$ for a sufficiently large
$\gamma$, a suitable energy $z(t)$ and $C$ independent of the
regularization parameter $\eps$ (in the proof of Lemma~\ref{transfert},
we need to establish such an inequality with  $\gamma=1$, see \eqref{gronwall} below). 
However the situation is the context of  \eqref{eq:regMKPII_mod} is technically more
complicated because of the presence of the parabolic regularization terms.
In the reasoning one uses the a priori estimates obtained for the equation satisfied by the Miura transform of $\tilde{Q}_c+w$, where $w$ is a
solution of \eqref{eq:regMKPII_mod}. 
It turns out that if $w$ solves \eqref{eq:regMKPII_mod} then $u(t,x,y)=(M^c_{+}(\tilde{Q}_c+w))(t,x-3ct,y)$ solves
\begin{equation}\label{eq:regMKPII_mod_miura}
\begin{split}
&  \pd_t u+\eps\pd_x^4u+\eps^5\pd_y^4u+\pd_x^3u +3\pd_x^{-1}\pd_y^2 u+3\partial_{x}(u^2)=
\\ &
-4\eps[\pd_x^2(\pd_x v)^2-\frac{1}{2}(\pd_x^2 v)^2]
-4\eps^5[\pd_y^2(\pd_y v)^2-\frac{1}{2}(\pd_y^2 v)^2]
-2\eps(\pd^4_x \hat{Q}_c)v+\eps\pd_x^5\hat{Q}_c\,, 
  \end{split}
\end{equation}
where $v=\hat{Q}_c+w(t,x-3ct,y)$ and $\hat{Q}_{c}=Q(x-2ct)$.
The equation \eqref{eq:regMKPII_mod_miura} is of KP type (with nicer nonlinearity compared to \eqref{eq:regMKPII_mod})
and thus the traditional energy estimates for the KP equations (or first order hyperbolic PDE's) apply to it.
These estimates together with a reasoning in the spirit of Lemma~\ref{transfert} below transform the energy estimates for
\eqref{eq:regMKPII_mod_miura}  to energy estimates for \eqref{eq:regMKPII_mod}.
The additional term in \eqref{eq:regMKPII_mod_miura} with respect to \cite{KMa} is $-2\eps(\pd^4_x \hat{Q}_c)v+\eps\pd_x^5\hat{Q}_c$
which is of lower order compared to the other terms in the right hand-side of \eqref{eq:regMKPII_mod_miura} and thus the analysis 
performed in \cite{KMa} is not affected. The a priori uniform in $\eps$ estimates for \eqref{eq:regMKPII_mod} imply the 
local well-posedness for \eqref{MKPII_mod} by a classical compactness argument.
 
Once local solutions of \eqref{MKPII_mod} are obtained as in \cite{KMa},
one has a global solution of \eqref{MKPII_mod}-\eqref{data_mod} thanks to
the global well-posedness of KP-II posed on $\R_x\times\T_y$  proved in \cite{MST} and the following lemma.
\begin{lemma}\label{transfert}
Suppose that $w\in Z$ is a solution of \eqref{MKPII_mod}-\eqref{data_mod}
on a time interval $[0,T)$. Let $u=M^c_{+}(\tilde{Q}_c+w)$
and suppose that 
$$
\sup_{t\in [0,T)}\|u(t,\cdot)\|_{H^8(\R_x\times\T_y)}<\infty.
$$
Then 
$$
\sup_{t\in [0,T)}\|w(t,\cdot)\|_{Z}<\infty\, .
$$
\end{lemma}
\begin{proof}[Proof of Lemma~\ref{transfert}]
Here we will use the method of Lemma~9 in \cite{KMa} by incorporating
a small modification coming from the presence of $c$ in the Miura transform
$M_{+}^c$. Since $M_+^c(Q_c)=\varphi_c$, we have
$$
M_{+}^c(\tilde{Q}_c+w)=\varphi_c(x+ct)+\partial_{x}w+\partial_{x}^{-1}\pd_yw-w^2-2\tilde{Q}_cw\,.
$$
Thus
\begin{equation}\label{fil1}
\sup_{t\in[0,T)}
\|\partial_{x}w+\partial_{x}^{-1}\partial_y w-w^2-2\tilde{Q}_c w\|_{L^2}\leq C\,.
\end{equation}
Combining the fact that $(w_{x},\partial_{x}^{-1}\partial_y w)=0,\quad  (w_x,w^2)=0$
and
$$
-2(w_x,\tilde{Q}_cw)=\int_{\R_x\times\T_y}Q_{c}'(x+ct)w^2(x,y)dxdy>0,
$$
with \eqref{fil1}, we have
\begin{equation}\label{vajna1}
\sup_{t\in[0,T)}\|w_x\|_{L^2}
+
\sup_{t\in[0,T)}\|\partial_{x}^{-1}\partial_y w-w^2-2\tilde{Q}_c w\|_{L^2}\leq C\,.
\end{equation}
Using Lemma~\ref{Besov} and the bound for $\|w_x\|_{L^2}$ we have just obtained, we have for $t\in [0,T)$,
$$
\|\partial_{x}^{-1}\partial_y w\|_{L^2}\leq C(\|w\|_{L^4}^2+\|w\|_{L^2})
\leq C(\|w\|_{L^2}^{\frac{1}{2}}
\|\partial_{x}^{-1}\partial_y w\|_{L^2}^{\frac{1}{2}}
+\|w\|_{L^2})
$$
which in turn implies that for $t\in [0,T)$,
\begin{equation}\label{vajna2}
\|\partial_{x}^{-1}\partial_y w\|_{L^2}\leq C\|w\|_{L^2}\,.
\end{equation}
We now obtain estimates for $\|w\|_{L^2}$. We multiply \eqref{MKPII_mod} by $w$ and integrate over
$\R\times\T$ to have after some integrations by parts
$$
\frac{1}{2}\frac{d}{dt}\|w\|_{L^2}^2
=-6\int_{\R_x\times\T_y}\tilde{Q}_c'\, w\,\partial_{x}^{-1}\partial_yw
+6\int_{\R_x\times\T_y}\tilde{Q}_c \tilde{Q}_c' w^2
+2\int_{\R_x\times\T_y}\tilde{Q}_c' w^3.
$$
Using Lemma~\ref{Besov}, \eqref{vajna1} and \eqref{vajna2}, we have
$$
\Big|\int_{\R_x\times\T_y}\tilde{Q}_c' w^3\Big|\leq \|\tilde{Q}_c'\|_{L^2(\R_x\times\T_y)}\|w\|_{L^6}^3
\leq C\|\partial_{x}^{-1}\partial_y w\|_{L^2}
\leq C\|w\|_{L^2}\,.
$$
Using the last estimate and \eqref{vajna2}, we have
\begin{equation}\label{gronwall}
\frac{d}{dt}\|w\|^2_{L^2} \leq C(\|w\|_{L^2}+1)^2\,.
\end{equation}
Therefore thanks to Gronwall's lemma, we have
$$\sup_{t\in[0,T)}\|w\|_{L^2}\leq C.$$
Therefore, we have obtained the needed bounds for 
$\|w\|_{L^2}$, $\|\partial_xw\|_{L^2}$ and $\|\partial_{x}^{-1}\partial_y w\|_{L^2}$, i.e.
$
E_0(w)\leq C,
$
where 
$$
E_0(w)=\|w\|_{L^2}+\|\partial_xw\|_{L^2}+\|\partial_{x}^{-1}\partial_y w\|_{L^2}\,.
$$
We next estimate higher derivatives. Write
\begin{eqnarray*}
\partial_{x}u & = &\partial_x\varphi_c(x+ct)+\partial^2_{x}w+\pd_yw-2w\partial_x w-2\partial_x(\tilde{Q}_cw),
\\
\partial_{y}u & = & \partial_{x}\partial_{y}w+\partial_{x}^{-1}\partial_{y}^{2}w-2w\partial_{y}w -2\tilde{Q}_c\partial_{y}w\,.
\end{eqnarray*}
Set
$$
E_1(w)=\|w\|_{L^2}+\|\partial_{x}^2 w\|_{L^2}+
\| \partial_{x}\partial_{y}w\|_{L^2}+\|\partial_{x}^{-1}\partial^2_y w\|_{L^2}\,.
$$
Using the orthogonality between $\partial^2_{x}w$ and $\pd_yw$ and also between $\partial_{x}\partial_{y}w$ and
$\partial_{x}^{-1}\partial_{y}^{2}w$, we obtain that
$$
E_{1}(w)\leq C(1+\|w\partial_x w\|_{L^2}+\|w\partial_y w\|_{L^2}+\|\partial_x w\|_{L^2}+\|\partial_y w\|_{L^2}).
$$
Next, thanks to an elementary interpolation inequality, we get
$$
\|\partial_x w\|_{L^2}+\|\partial_y w\|_{L^2}\leq C(E_{0}(w))^{\frac{1}{2}}(E_1(w))^{\frac{1}{2}}\leq C(E_{1}(w))^{\frac{1}{2}}\,.
$$
Next, we write by invoking Lemma~\ref{Besov},
\begin{eqnarray*}
\|w\partial_x w\|_{L^2} & \leq &\|w\|_{L^6}\|\partial_x w\|_{L^3}\leq
CE_{0}(w)\|\partial_x w\|_{L^2}^{\frac{1}{2}}\|\partial_{x}^2 w\|_{L^2}^{\frac{1}{3}}\|\partial_{y} w\|_{L^2}^{\frac{1}{6}}
\\
& \leq &
C(E_0(w))^{\frac{5}{4}}(E_1(w))^{\frac{3}{4}}\leq C(E_1(w))^{\frac{3}{4}}\,.
\end{eqnarray*}
Similarly
\begin{eqnarray*}
\|w\partial_y w\|_{L^2} & \leq &\|w\|_{L^6}\|\partial_y w\|_{L^3}\leq
CE_{0}(w)\|\partial_y w\|_{L^2}^{\frac{1}{2}}\|\partial_{x}\partial_{y} w\|_{L^2}^{\frac{1}{3}}\|\partial_{x}^{-1}\partial^2_{y} w\|_{L^2}^{\frac{1}{6}}
\\
& \leq &
C(E_0(w))^{\frac{5}{4}}(E_1(w))^{\frac{3}{4}}\leq C(E_1(w))^{\frac{3}{4}}\,.
\end{eqnarray*}
In summary, we get $E_{1}(w)\leq C(1+(E_1(w))^{\frac{3}{4}})$ which gives
$E_1(w)\leq C$. Observe that $\|w\|_{L^\infty}\leq C E_{1}(w)$, i.e. we already have a control on the $L^\infty$ norm. 

Next, we write
\begin{eqnarray*}
\partial_{x}^8u & = &\partial^8_x\varphi_c(x+ct)+\partial^9_{x}w+\partial_{x}^7\pd_yw-\partial^8_x(w^2)-2\partial^8_x(\tilde{Q}_cw),
\\
\partial^8_{y}u & = & \partial_{x}\partial^8_{y}w+\partial_{x}^{-1}\partial^9_{y}w-\partial^8_{y}(w^2) -2\tilde{Q}_c\partial^8_{y}w\,.
\end{eqnarray*}
Set
$$
E_8(w)=\|w\|_{L^2}+\|\partial_{x}^{-1}\partial_y w\|_{L^2}+
\|\partial_{x}^9 w\|_{L^2}+\|\partial_x^7 \partial_{y} w\|_{L^2}+
\| \partial_{x}\partial^8_{y}w\|_{L^2}+\|\partial_{x}^{-1}\partial^9_y w\|_{L^2}\,.
$$
By invoking an orthogonality argument and the bounds we have already obtained, we can write
$$
E_{8}(w)\leq C(1+\|\partial^8_x (w^2)\|_{L^2}+\|\partial^8_y (w^2)\|_{L^2}+\|\partial^8_x w\|_{L^2}+\|\partial^8_y w\|_{L^2}).
$$
By invoking  the elementary inequality
\begin{equation}\label{elem}
a^{\theta}\,b^{1-\theta}\leq a+b,\qquad \forall\, a\geq0,\,\,\forall\, b\geq0,\,\, \forall\, \theta\in (0,1)
\end{equation}
in conjugation with the Fourier transform and the controls we have already obtained, we obtain that there exists
$\theta\in (0,1)$ such that. 
$$
\|\partial^8_x w\|_{L^2}+\|\partial^8_y w\|_{L^2}\leq C(E_{8}(w))^{\theta}\,.
$$
Next, we use a classical multiplicative inequality to get there exists $\theta\in (0,1)$ such that
\begin{eqnarray*}
\|\partial^8_x (w^2)\|_{L^2}+\|\partial^8_y (w^2)\|_{L^2} & \leq  &
C\|w^2\|_{H^8}\leq C\|w\|_{L^\infty}\|w\|_{H^8}
\\
& \leq & C(1+\|\partial_x^8 w\|_{L^2}+\|\partial_y^8 w\|_{L^2})\leq
C(1+(E_8(w))^\theta).
\end{eqnarray*}
Therefore, we obtain that $E_{8}(w)\leq C(1+(E_8(w))^\theta)$ for some $\theta\in (0,1)$. This in turns implies that $E_8(w)\leq C$, by a suitable use 
of \eqref{elem}.
We finally observe that $\|w\|_{Z}\leq C E_8(w)$. This completes the proof of Lemma~\ref{transfert}.
\end{proof}
Observe that using \eqref{algebra}, we infer that if $w\in Z$, then 
$$
M^{\pm}_c(Q_c+w)\in H^8(\R_x\times\T_y)\,.
$$
Once global solutions are established,
the conservation laws are obtained due to the following lemma.
\begin{lemma}\label{luc}
Let $u\in C(\R;H^8(\R_x\times\T_y))$ be a solution of the KP-II equation
\eqref{eq:KPII}. Then
$$
\|u(t,\cdot)\|_{L^2(\R_x\times\T_y)}=\|u(0,\cdot)\|_{L^2(\R_x\times\T_y)},\quad \forall\, t\in\R.
$$
\end{lemma}
The proof Lemma~\ref{luc} can be obtained by an argument due to Molinet (see \cite{Molinet}).
A similar argument may also be found in \cite[p.~785]{MST07}.
This completes the proof of Proposition~\ref{prop:wpmkp}.
\end{proof}
\end{proposition}
\section{The Miura transform $M^c_{+}$ in a neighborhood of $Q_c$}
\label{sec:miura}
It turns out that the Miura transform $M^c_{+}$ defines a bijection between a neighborhood of $(c,Q_c)$ and a neighborhood of $\varphi_c$.
\begin{proposition}\label{prop:bijMiura}
For every $\varepsilon>0$, there exists a $\delta>0$ such that if
$\|u\|_{L^2}<\delta$, there exists a unique
$(k,v)\in\R\times\mathcal{E}(\R_x\times\T_y)$ satisfying
\begin{equation}\label{eq:Miura2}
|k-c|<\varepsilon, \quad\|v\|_{\mathcal{E}(\R_x\times\T_y)}<\varepsilon,\quad
M^k_{+}(Q_k+v)=\varphi_c+u.
\end{equation}
Moreover, the map $L^2(\R_x\times\T_y)\ni u\mapsto (k,v)\in\R\times\mathcal{E}(\R_x\times\T_y)$
is of class $C^1$.
\end{proposition}
To prove Proposition~\ref{prop:bijMiura}, we need to investigate
a linearized operator of the Miura transform $M_+^c$ around $Q_c$.
Let $(\varphi_c)^{\perp}$ be a subspace of $L^2(\R_x\times\T_y)$ defined by
$$(\varphi_c)^{\perp}\equiv (u\in L^2(\R_x\times\T_y)\,:\, 
(u,\varphi_c)_{L^2(\R_x\times\T_y)}=0)\,.$$
We will show that the Fr\'echet derivative 
$\nabla_u M_+^c(Q_c):\mathcal{E}\to(\varphi_c)^{\perp}$ is bicontinuous
and bounded from below.
\begin{lemma}
\label{lem:heatop}
Let $c>0$ and consider $L_c\equiv\pd_x+\pd_x^{-1}\pd_y-2Q_c(x)$  and its formal adjoint
$\mathcal{L}_c\equiv-\pd_x+\pd_x^{-1}\pd_y-2Q_c(x)$  as 
bounded operators  from $\mathcal{E}(\R_x\times\T_y)$ to $L^2(\R_x\times\T_y)$.
Then $\ker(L_c)=\{0\}$ and $\ker(\mathcal{L}_c)=\spann\{\varphi_c\}$.
Moreover $L_c$ is a Fredholm operator and
$\operatorname{Range}(L_c)=(\varphi_c)^\perp$.
In addition,
\begin{equation}\label{fred}
\|w\|_{\mathcal{E}(\R_x\times\T_y)}\leq C\|L_cw\|_{L^2(\R_x\times\T_y)},
\end{equation}
where $C$ is a positive constant that does not depend on $w$.
\end{lemma}
\begin{proof}[Proof of Lemma~\ref{lem:heatop}]
First, we remark that if $u\in \ker(L_c)$ or $u\in \ker(\mathcal{L}_c)$,  then $\pd_x^i\pd_y^ju\in L^2(\R_x\times\T_y)$ for every $i\ge 0$ and $j\ge0$
thanks to an elliptic regularity argument.
We will give the proof only for $c=2$, the case of a general $c$ being
the same modulo some direct modifications.
Let $Q\equiv Q_2$, $\varphi\equiv \varphi_2$, $L\equiv L_2$ and $\mathcal{L}
=\mathcal{L}_2$.
\par
Suppose $u\in \ker(L)$. Then $u$ satisfies
\begin{equation}
  \label{eq:-heat}
u_y+u_{xx}=2(Qu)_x.  
\end{equation}
Using  the last equation, we obtain that $u$ satisfies
$$
\frac{1}{2}\frac{d}{dy}\int_{\R}u^2(x,y)dx=\int_\R(u_x^2(x,y)+Q'(x)u^2(x,y))dx\,.
$$
We next integrate the above identity over $\T$ to have
$$
0=\int_{\T}\int_\R(u_x^2(x,y)+Q'(x)u^2(x,y))dxdy\,.
$$
Combining the above with the fact that $Q'(x)>0$ for every $x\in\R$, we have
$u=0$. Thus we obtain that $\ker(L)=\{0\}$.
\par
The study of $\ker(\mathcal{L})$ is more intricate. 
Suppose $u\in\ker(\mathcal{L})$. Then $u$ is a solution to a heat equation
\begin{equation}
  \label{eq:heat}
  u_y=(u_x+2Qu)_x,
\end{equation}
and $2\pi$-periodic in $y$. A direct computation shows that \eqref{eq:heat}
has $y$-independent solutions $\{\alpha Q'(x)\,|\, \alpha \in\R\}$.
We will show that \eqref{eq:heat} has no other solution which is periodic in the $y$-variable.
Let 
$$
V(y)=\int_\R\Big(\frac12u_x^2(x,y)-(Q'(x)-2Q^2(x))u^2(x,y)\Big)dx.
$$
If $u\in L^2(\R\times \T)$ is a smooth (in the Sobolev scale) solution of \eqref{eq:heat},
\begin{equation*}
  \begin{split}
V'(y)=& -\int_\R u_y(x,y)\big(u_{xx}(x,y)+2(Q'(x)-2Q^2(x))u(x,y)\big)dx
\\=&-\int_\R u_y(x,y)\big(u_y(x,y)-2Q(x)(u_x(x,y)+2Q(x)u(x,y))\big)dx
\\=&-\int_\R \big(u_y^2(x,y)+Q'(x)(u_x(x,y)+2Q(x)u(x,y))^2\big)dx\leq 0.     
  \end{split}
\end{equation*}
Integrating the last identity over $\T$, we have $u_y=u_x+2Qu=0$.
Thus $u$ is independent of $y$ and by solving the ODE $u_x+2Qu=0$,
we obtain that $\ker(\mathcal{L})=\spann\{Q'\}$.
\par
Finally, we will show $R(L)=L^2(\R_x\times\T_y)\cap (Q')^\perp$
to prove that $L$ is Fredholm.
Since $L$ is formally an adjoint operator of $\mathcal{L}$ and
$\ker(\mathcal{L})=\spann\{Q'\}$, we have
$R(L)\subset L^2(\R_x\times\T_y)\cap (Q')^\perp$.
Thus it suffices to show that $Lu=f$ has a solution
$u\in\mathcal{E}(\R_x\times\T_y)$ for any $f\in L^2(\R\times\T)\cap (Q')^\perp$.

Suppose that $u$ satisfies $Lu=f$ with $f\in L^2(\R_x\times\T_y)\cap (Q')^\perp$.
Let us expand $f$ and $u$ into Fourier series in $y$:
$$
f(x,y)=\frac{1}{\sqrt{2\pi}}\sum_{n\in\Z}f_n(x)e^{iny},\quad
u(x,y)=\frac{1}{\sqrt{2\pi}}\sum_{n\in\Z}u_n(x)e^{iny}\,.
$$
Then we find that $u_n$ and $f_n$ satisfy the equations
\begin{equation}
  \label{eq:fouriermode}
  \begin{split}
&  u_0'-2Qu_0=f_0\,,\quad \int_{-\infty}^\infty f_0(x)Q'(x) dx=0.\\
&  u_n''+inu_n-2(Qu_n)'=f_n' \quad\text{if $n\ne0$},
  \end{split}
\end{equation}
and $f_{n}\in L^2(\R)$ for every $n\in\Z$.
The first equation in \eqref{eq:fouriermode} can be solved explicitly as
\begin{equation}\label{u0_formula}
u_0(x)=-\int_x^\infty\frac{\varphi(t)f_0(t)}{\varphi(x)}dt
=\int_{-\infty}^x\frac{\varphi(t)f_0(t)}{\varphi(x)}dt\,.
\end{equation}
Note that $\varphi=2Q'$. We have the following bound on $u_0$.
\begin{claim}\label{u0}
There exists $C$ such that for every $f_0\in L^2(\R)$,
$$
\|u_0\|_{H^1(\R)}\leq C\|f_0\|_{L^2}\,.
$$
\end{claim}
\begin{proof}
Since  $\varphi(t)/\varphi(x)\le 4e^{-2|x-t|}$ for $|t|\ge |x|$,
it follows that
$$|u_0(x)|\le 4 \int_\R e^{-2|x-t|}|f_0(t)|dt\,.$$
Applying Young's inequality to the above, we obtain
$\|u_0\|_{L^2}\le C\|f_0\|_{L^2}.$
The estimate for $\|u'_0\|_{L^2}$ follows from the equation satisfied by $u_0$ and the previous estimate.
This completes the proof of Claim~\ref{u0}.
\end{proof}
Next, we will show solvability of the second equation in \eqref{eq:fouriermode}.
The key point is to show that the homogeneous problem
\begin{equation}\label{pb_hom}
u''+inu-2(Qu)'=0
\end{equation}
has no nontrivial spatially localized solution. Indeed, we have the following statement.
\begin{claim}\label{empty}
Let $n\in\Z$. Then \eqref{pb_hom} has no nontrivial solution $u$
which belongs to $L^2(\R)$.
\end{claim}
\begin{proof}
Let $u\in L^2(\R)$ be a solution of \eqref{pb_hom}.
Then by an elliptic regularity argument, we have $u\in H^2(\R)$. 
Moreover, we see that $e^{iny}u(x)\in \mathcal{E}(\R_x\times\T_y)$ and
$L_c(e^{iny}u(x))=0$. Since $\ker(L_c)=\{0\}$,
it follows that $u=0$. This completes the proof of Claim~\ref{empty}.
\end{proof}
By standard ODE arguments, 
since $\lim_{x\to\pm\infty}Q'(x)=0$ and $\lim_{x\to\pm\infty}Q(x)=\pm1$,
the equation \eqref{pb_hom} has fundamental systems $\{w_n^+(x),v_n^+(x)\}$
and $\{w_n^-(x),v_n^-(x)\}$ such that
\begin{align*}
& w_n^+(x)\sim e^{\mu_n^+x},\quad v_n^+(x)\sim e^{\lambda_n^+x}
\quad\text{as $x\to\infty$,}\\
& w_n^-(x)\sim e^{\mu_n^-x},\quad v_n^-(x)\sim e^{\lambda_n^-x}
\quad\text{as $x\to-\infty$,}
\end{align*}
where $\mu_n^\pm=\pm1+\sqrt{1-in}$ and $\lambda_n^\pm=\pm1-\sqrt{1-in}$.
Here and in the sequel the square roots are taken so that the real part of the result is non negative.
Since ${\rm Re} (\mu_n^\pm)>0>{\rm Re}(\lambda_n^\pm)$ for every $n\ne0$,
we see that $v_n^+(x)$ decays exponentially as $x\to\infty$ 
and $w_n^-(x)$ decays exponentially as $x\to-\infty$,
whereas that $w_n^+(x)$ grows exponentially as $x\to\infty$ 
and $v_n^-(x)$ grows exponentially as $x\to-\infty$.
Using Claim~\ref{empty}, we see that $v_n^+(x)$ and $w_n^-(x)$
are linearly independent for $n\ne0$. In other words,
$v_n^+(x)=O(v_n^-(x))$ as $x\to-\infty$ and $w_n^-(x)=O(w_n^+(x))$ as
$x\to\infty$ if $n\ne0$. The Green kernel is given by
\begin{gather}
  \label{eq:green}
G_n(x,t)=\left\{\begin{aligned}
& -\frac{v_n^+(x)w_n^-(t)}{W(v_n^+,w_n^-)(t)}\quad\text{for $x>t$,}\\
& -\frac{w_n^-(x)v_n^+(t)}{W(v_n^+,w_n^-)(t)}\quad\text{for $x<t$},\\    
\end{aligned}\right.
\end{gather}
where
$$W(v_n^+,w_n^-)(x)=\begin{vmatrix}
v_n^+(x) & w_n^-(x)\\ \pd_xv_n^+(x) & \pd_xw_n^-(x)
\end{vmatrix}=\cosh^2(x)W(v_n^+,w_n^-)(0).
$$
Observe that thanks to the above properties of the Green kernel, 
the kernel $G_n$ enjoys a pointwise bound
$$
|G_n(x,t)|\lesssim \exp(-({\rm Re}\sqrt{1-in}-1)|x-t|)\,.
$$
Thus for every $n\ne0$, the second equation of \eqref{eq:fouriermode} has a
solution given by
$$
u_n(x)=\int_{\R}G_{n}(x,t)f_n'(t)dt\,.
$$
We now observe that 
$$\pd_tG_n(x,t)=\left\{\begin{aligned}
& 2Q(t)G_n(x,t)-\frac{v_n^+(x)\pd_tw_n^-(t)}{W(v_n^+,w_n^-)(t)}
\quad\text{for $x>t$,}\\
& 2Q(t)G_n(x,t)-\frac{w_n^-(x)\pd_tv_n^+(t)}{W(v_n^+,w_n^-)(t)}
\quad\text{for $x<t$},\\    
\end{aligned}\right.$$
and 
\begin{equation*}
|\pd_tG_n(x,t)|\lesssim \exp(-({\rm Re}\sqrt{1-in}-1)|x-t|)\,.
\end{equation*}
Therefore using integration by parts
(and an approximation argument for $f_n$ by $C_0^\infty$ functions), we get
$$
u_n(x)=-\int_{\R}\partial_{t}G_{n}(x,t)f_n(t)dt\,.
$$
Differentiating the above equation, we have
$$\pd_xu_n(x)=f_n(x)+\int_\R G_n^1(x,t)f_n(t)dt,$$
where 
$$G_n^1(x,t)=
\left\{\begin{aligned}
& 2Q(t)\pd_xG_n(x,t)-\frac{\pd_xv_n^+(x)\pd_tw_n^-(t)}{W(v_n^+,w_n^-)(t)}
\quad\text{for $x>t$,}\\
& 2Q(t)\pd_xG_n(x,t)-\frac{\pd_xw_n^-(x)\pd_tv_n^+(t)}{W(v_n^+,w_n^-)(t)}
\quad\text{for $x<t$}.\\    
\end{aligned}\right.$$
Obviously, 
\begin{equation*}
|G_n^1(x,t)|\lesssim \exp(-({\rm Re}\sqrt{1-in}-1)|x-t|)\,.
\end{equation*}
Thus we obtain that $u_n\in H^1$ and
\begin{equation}\label{energia}
\|u_n\|_{H^1}+|n|\|\partial_{x}^{-1}u_n\|_{L^2}\leq C(n)\|f_{n}\|_{L^2}
\end{equation}
follows from the Young's inequality and the relation
$in\pd_x^{-1}u_n=2Qu_n-\pd_xu_n$.
\par

Now we will show that $C(n)$ can be chosen uniformly in $n$.
Let $T_n(u)\equiv 2\pd_x(in+\pd_x^2)^{-1}(Qu)$ and
$g_n\equiv\pd_x(in+\pd_x^2)^{-1}f_n$.
Then the second equation of \eqref{eq:fouriermode} can be rewritten as
$$
u_n=T_n(u_n)+g_n.
$$
Since $\|T_n\|_{B(H^1)}=O(1/\sqrt{|n|})$ by the Plancherel theorem,
there exists an $n_0$ such that $\|T_n\|_{B(H^1)}\le 1/2$ for every $|n|\ge n_0$. Hence there exists a positive
number $C$ such that for every $|n|\ge n_0$ 
$$\|u_n\|_{H^1}\le C\|g_n\|_{H^1}\le C\|f_n\|_{L^2}.$$
Furthermore, there exists a $C'>0$ such that for every $|n|\ge n_0$,
\begin{align*}
\|\pd_x^{-1}(T_nu_n+g_n)\|_{L^2} \le & 
2\|(\pd_x^2+in)^{-1}(Qu_n)\|_{L^2}+\|(\pd_x^2+in)^{-1}f_n\|_{L^2}
\\ \le & \frac{1}{|n|}(2\|u_n\|_{L^2}+\|f_n\|_{L^2})\le \frac{C'\|f_n\|_{L^2}}{|n|},
\end{align*}
whence $\|\pd_x^{-1}u_n\|\le C'\|f_n\|_{L^2}/|n|$ for $|n|\ge n_0$.
Therefore the constant $C(n)$ involved in \eqref{energia} is uniform in $n$.
Combining \eqref{energia} and Claim~\ref{u0}, we have
\begin{align*}
\|u\|_{\mathcal{E}(\R_x\times\T_y)}^2=& \sum_{n\in\Z}
(\|u_n\|_{H^1(\R)}^2+n^2\|\pd_x^{-1}u_n\|_{L^2(\R)}^2)
\\ \leq & C\sum_{n\in\Z}\|f_n\|_{L^2(\R)}^2=C\|f\|_{L^2(\R_x\times\T_y)}^2.
\end{align*}
Thus we have proved that $u\in\mathcal{E}(\R_x\times\T_y)$ and
$L:\mathcal{E}(\R_x\times\T_y)\to (\varphi)^\perp$ is surjective.
Finally, the estimate \eqref{fred} follows readily from the open mapping
theorem. This completes the proof of Lemma~\ref{lem:heatop}.
\end{proof}
Now we are in position to prove Proposition \ref{prop:bijMiura}.
\begin{proof}[Proof of Proposition \ref{prop:bijMiura}]
Since $M^c_{+}(Q_c)=\varphi_c$, \eqref{eq:Miura2} can be rewritten as
\begin{equation}
\label{eq:Miura3}
  \varphi_c-\varphi_k=L_cv+2v(Q_c-Q_k)-v^2-u.
\end{equation}
Let $P$ be a projection from $L^2(\R_x\times\T_y)$ to its orthogonal subspace
$(\varphi)^\perp$ defined by
$$
Pf\equiv f-\frac{(f,Q_c')_{L^2}}{\|Q_c'\|_{L^2}^2}Q_c'\,.
$$
Next, we define $F_1$ and $F_2$ by
\begin{align*}
  & F_1(u,v,k)=(\varphi_k-\varphi_c+2v(Q_c-Q_k)-v^2-u,\varphi_c)_{L^2},\\
  & F_2(u,v,k)=L_cv+P(\varphi_k-\varphi_c+2v(Q_c-Q_k)-v^2-u).
\end{align*}
Since $\mathcal{L}_c\varphi_c=0$, \eqref{eq:Miura3} holds if and only if
$F_1(u,v,k)=0$ and $F_2(u,v,k)=0$.
We consider $(F_1(u,v,k),F_2(u,v,k))$ as a $C^1$ map from 
$L^2(\R_x\times\T_y)\times \mathcal{E}(\R_x\times\T_y)\times\R_+$ to
$\R\times (\varphi_c)^{\perp}$.
Observe that $F_1(0,0,c)=0$, $F_2(0,0,c)=0$. Next we compute
$$
{\mathcal A}\equiv
\begin{pmatrix}
\pd_kF_1 &  \pd_vF_1\\ \pd_kF_2 & \pd_vF_2\end{pmatrix}\biggm|_{(u,v,k)=(0,0,c)}
=\begin{pmatrix}
  \big((\partial_k\varphi_k)|_{k=c},\varphi_c\big) & 0\\ P((\pd_k\varphi_k)|_{k=c}) & L_c
\end{pmatrix}\,.
$$
Since
$$
\big((\partial_k\varphi_k)|_{k=c},\varphi_c\big)
=\frac{3}{4c}\|\varphi_c\|_{L^2}^2\neq 0,
$$
we obtain that ${\mathcal A}$ is a bicontinuous bijection from $\R\times \mathcal{E}(\R_x\times\T_y)$ to $\R\times (\varphi_c)^{\perp}$
by using Lemma~\ref{lem:heatop}. Therefore the assertion of Proposition~\ref{prop:bijMiura} follows from the implicit function theorem.
\end{proof} 

Next we will further investigate the linearized operator of the Miura transform $M_-^c$. 
\begin{lemma}
\label{lem:heatop2}
Let $c>0$ and $\mathcal{L}_c\equiv-\pd_x+\pd_x^{-1}\pd_y-2Q_c(x)$ be considered as
a bounded operator from $\mathcal{E}(\R_x\times\T_y)$ to $L^2(\R_x\times\T_y)$.
Then $\mathcal{L}_c$ is Fredholm. More precisely,
$\ker(\mathcal{L}_c)=\spann\{\varphi_c\}$,  $\operatorname{Range}(\mathcal{L}_c)=L^2$ and 
$\mathcal{L}_c:\mathcal{E}\cap(Q_c')^\perp\to L^2$ has a bounded inverse.
\end{lemma}
To prove Lemma \ref{lem:heatop2}, we need the following.
\begin{claim}
  \label{cl:u0}
Let $I(f)(x)=\varphi(x)\int_0^xf(t)\varphi(t)^{-1}dt.$
Then there exists a positive constant $C$ such that
$$\|I(f)\|_{H^1(\R)}\le C\|f\|_{L^2(\R)}.$$
\end{claim}
\begin{proof}
Making use of
\begin{equation}
  \label{eq:vx/vt}
\varphi(x)/\varphi(t)\le 4e^{-2|x-t|},\quad
\left|\varphi'(x)/\varphi(t)\right|\le 8e^{-2|x-t|}
\quad\text{for $t\in[-|x|,|x|]$,}
\end{equation}
we have $\|I(f)\|_{H^1}\le C\|f\|_{L^2}$ in exactly the same way as Claim~\ref{u0}.  
\end{proof}
\begin{proof}[Proof of Lemma \ref{lem:heatop2}]
We will give the proof only for $c=2$ for the sake of simplicity.
Let $\mathcal{L}=\mathcal{L}_2$.
We have $\ker(\mathcal{L})=\spann\{Q'\}$ from Lemma \ref{lem:heatop}.
\par
Next, we will show that $\mathcal{L}u=f$ has a solution $u\in\mathcal{E}$
for any $f\in L^2$. Let us expand $u$ and $f$ into Fourier series as
\begin{equation*}
u(x,y)=\sum_{n\in\Z}u_n(x)e^{iny}\,,\quad f(x,y)=\sum_{n\in\Z}f_n(x)e^{iny}.  
\end{equation*}
Then $f_n\in L^2(\R)$ for every $n\in\Z$ and
\begin{align}
  \label{eq:Lf0}
&  -u_0'-2Qu_0=f_0,\\
\label{eq:Lfn}
&  -u_n''+inu_n-2(Qu_n)'=f_n' \quad\text{if $n\ne0$.}
\end{align}
If $u_0$ is a solution of \eqref{eq:Lf0}, then
$u_0(x)=\alpha\varphi(x)+I(f_0)(x)$ for an $\alpha\in\R$.
We remark that $\alpha$ is uniquely determined and
$L^2(\R)\ni f_0\mapsto \alpha\in\R$ is continuous
if we impose the orthogonality condition $\int u_0Q'dxdy=0$.
\par

Let us now observe that for $n\neq0$, the equation \eqref{eq:Lfn} has no nontrivial solution $u$
which belongs to $L^2(\R)$. Indeed, let $u\in L^2(\R)$ be a solution of \eqref{eq:Lfn}.
Then by elliptic regularity, we have $u\in H^2(\R)$.  Moreover, we see that $e^{iny}u(x)\in \mathcal{E}(\R_x\times\T_y)$ and
$\mathcal{L}(e^{iny}u(x))=0$.  Since $\ker(\mathcal{L})=\spann\{Q'\}$ it follows that $u=0$ unless $n=0$.
Next, we will solve \eqref{eq:Lfn}.
By standard ODE arguments, since $\lim_{x\to\pm\infty}Q'(x)=0$
and $\lim_{x\to\pm\infty}Q(x)=\pm1$, the equation 
\begin{equation}
  \label{eq:kern}
  -u''+inu-2(Qu)'=0
\end{equation}
has fundamental systems $\{\tilde{w}_n^+(x),\tilde{v}_n^+(x)\}$
and $\{\tilde{w}_n^-(x),\tilde{v}_n^-(x)\}$ such that
\begin{align*}
& \tilde{w}_n^+(x)\sim e^{\tilde{\mu}_n^+x},
\quad \tilde{v}_n^+(x)\sim e^{\tilde{\lambda}_n^+x}
\quad\text{as $x\to\infty$,}\\
& \tilde{w}_n^-(x)\sim e^{\tilde{\mu}_n^-x},
\quad \tilde{v}_n^-(x)\sim e^{\tilde{\lambda}_n^-x}
\quad\text{as $x\to-\infty$,}
\end{align*}
where $\tilde{\mu}_n^\pm=\mp1+\sqrt{1+in}$
and $\tilde{\lambda}_n^\pm=\mp1-\sqrt{1+in}$.
Since ${\rm Re} (\tilde{\mu}_n^\pm)>0>{\rm Re}(\tilde{\lambda}_n^\pm)$
for every $n\ne0$,
\begin{align*}
& \lim_{x\to\infty}\tilde{w}_n^+(x)=\infty,
\quad \lim_{x\to\infty}\tilde{v}_n^+(x)=0,\\
& \lim_{x\to-\infty}\tilde{w}_n^-(x)=0,\quad
\lim_{x\to-\infty}\tilde{v}_n^-(x)=\infty.
\end{align*}
As in the proof of Claim \ref{empty},
we see that $\tilde{v}_n^+(x)$ and $\tilde{w}_n^-(x)$
are linearly independent for every $n\ne0$.
\par
Thus for $n\ne0$, the Green kernel of $-u''+inu-2(Qu)'$ is given by
\begin{gather}
  \label{eq:green2}
\widetilde{G}_n(x,t)=\left\{\begin{aligned}
&- \frac{\tilde{v}_n^+(x)\tilde{w}_n^-(t)}
{W(\tilde{v}_n^+,\tilde{w}_n^-)(t)}\quad\text{for $x>t$,}\\
& -\frac{\tilde{w}_n^-(x)\tilde{v}_n^+(t)}{W(\tilde{v}_n^+,\tilde{w}_n^-)(t)}
\quad\text{for $x<t$},\\\end{aligned}\right.
\end{gather}
where
$$W(\tilde{v}_n^+,\tilde{w}_n^-)(x)=\begin{vmatrix}
\tilde{v}_n^+(x) & \tilde{w}_n^-(x)\\
 \pd_x\tilde{v}_n^+(x) & \pd_x\tilde{w}_n^-(x)
\end{vmatrix}=\sech^2(x)W(\tilde{v}_n^+,\tilde{w}_n^-)(0).
$$
For every $n\ne0$,
\begin{equation}
  \label{eq:grest}
|\widetilde{G}_n(x,t)|+|\pd_t\widetilde{G}_n(x,t)|
\lesssim \exp(-{\rm Re}(\sqrt{1+in}-1)|x-t|).
\end{equation}
Repeating the arguments of the proof of Lemma~\ref{lem:heatop},
we obtain
\begin{equation}\label{eq:energia2}
\|u_n\|_{H^1}+|n|\|\partial_{x}^{-1}u_n\|_{L^2}\leq C(n)\|f_{n}\|_{L^2}
\end{equation}
and that the constant $C(n)$ involved in \eqref{eq:energia2} can be chosen
uniformly in $n$. Therefore
\begin{equation*}
\|u\|_{\mathcal{E}(\R_x\times\T_y)}^2= \sum_{n\in\Z}(\|u_n\|_{H^1(\R)}^2+n^2\|\pd_x^{-1}u_n\|_{L^2(\R)}^2)\leq C\|f\|_{L^2(\R_x\times\T_y)}^2,
\end{equation*}
where $C$ is a constant independent of $f\in L^2(\R_x\times\T_y)$.
Thus we have proved that $\mathcal{L}:\mathcal{E}(\R_x\times\T_y)\to
L^2(\R_x\times\T_y)$ is surjective.

Since $\ker(\mathcal{L})=\spann\{Q'\}$ and $\operatorname{Range}(\mathcal{L})
=L^2$, it follows from the open mapping theorem that
$\mathcal{L}:\mathcal{E}\cap(Q')^\perp\to L^2$ has a bounded inverse.
\end{proof}
Finally, we will investigate a property of $\mathcal{L}_c$ in a weighted space.  Let 
$$
\chi_\eps(x)\equiv\frac{1+\tanh (\eps x)}{2}
$$ 
and let $L^2_{\eps,x_0}$ and $\mathcal{E}_{\eps,x_0}$ be Banach spaces equipped with norms
\begin{gather*}
\|u\|_{L^2_{\eps,x_0}}=\Big(\int_{\R_x\times\T_y}\chi_{\eps}(x+x_0)u^2(x,y)dxdy\Big)^{\frac12},
\\ \|u\|_{\mathcal{E}_{\eps,x_0}}=\|u\|_{L^2_{\eps,x_0}}+\|\pd_xu\|_{L^2_{\eps,x_0}}
+\|\pd_x^{-1}\pd_yu\|_{L^2_{\eps,x_0}},
\end{gather*}
respectively. Roughly speaking, we will show that
$\mathcal{L}_c:\mathcal{E}_{\eps,x_0}\cap (Q_c')^\perp\to L^2_{\eps,x_0}$
is bounded from below for small $\eps>0$.
\begin{lemma}
\label{lem:heatop3}
Let $c>0$.  Then there exist positive constants $\eps_0$ and $C$ such that
for every $\eps\in(0,\eps_0)$, $x_0\in\R$ and 
$w\in \mathcal{E}(\R_x\times\T_y)\cap(Q_c')^\perp$,
\begin{equation}\label{eq:fred3}
\|\mathcal{L}_c w\|_{L^2_{\eps,x_0}(\R_x\times\T_y)} \ge C \|w\|_{\mathcal{E}_{\eps,x_0}}.
\end{equation}
\end{lemma}
\begin{proof}
Let $c=2$ for the sake of simplicity.
Suppose 
$
f(x,y)=\sum_{n\in\Z}f_n(x)e^{iny}\in L^2(\R_x\times\T_y)
$ 
and that
$u(x,y)=\sum_{n\in\Z}u_n(x)e^{iny}\in \mathcal{E}(\R_x\times\T_y)$
is a solution of $\mathcal{L}_cu=f$ satisfying $\int uQ'dxdy=0$.
Then 
\begin{gather*}
u_0(x)=\alpha\varphi(x)+I(f_0)(x),\\
u_n(x)=\int_\R \widetilde{G}_n(x,t)f'_n(t)dt\quad\text{for $n\ne0$,}
\end{gather*}
where $\alpha$ is a constant satisfying
$\alpha=-(I(f_0),\varphi)_{L^2}/\|\varphi\|_{L^2}^2$.
\par
By the definition of $\chi_\eps$,
\begin{equation}\label{eq:chch-1}
 \frac{\chi_\eps(x+x_0)}{\chi_\eps(t+x_0)}\leq 1+e^{2\eps(x-t)}\leq 2e^{2\eps|x-t|}\,.
\end{equation}
Combining \eqref{eq:vx/vt} and \eqref{eq:chch-1}, 
we have for $\eps\in(0,1)$,
$$\|\chi_{\eps}(x+x_0)I(f_0)\|_{L^2}
+\|\chi_{\eps}(x+x_0)\pd_xI(f_0)\|_{L^2}
\le C\|\chi_\eps(x+x_0)f_0\|_{L^2},$$
where $C$ is a positive constant depending only on $\eps$.
Hence there exists a $C>0$ such that for every $f\in L^2(\R)$ and $x_0\in\R$,
$$\|\chi_\eps(x+x_0)u_0\|_{L^2(\R)}+\|\chi_\eps(x+x_0)\pd_xu_0\|_{L^2(\R)}
\le C\|\chi_\eps(x+x_0)f_0\|_{L^2(\R)}.
$$
In view of \eqref{eq:grest} and \eqref{eq:chch-1}, we have
$$
|\chi_\eps(x+x_0)\pd_t^k\widetilde{G}_n(x,t)\chi_\eps(t+x_0)^{-1}|
\lesssim \exp((-{\rm Re}(\sqrt{1+in}-1)+2\eps)|x-t|).
$$
for $k=0$, $1$. As in the proof of Lemma~\ref{lem:heatop} after an analysis of
$\partial_x (\chi_\eps(x+x_0)\partial_t\widetilde{G}_n(x,t)\chi_\eps(t+x_0)^{-1})$,
we have that for $0<2\eps<{\rm Re}\sqrt{1+in}-1$,
there exist positive constants $C(n,\eps)$ such that if $n\ne0$,
$$
\|\chi_\eps(x+x_0)u_n\|_{L^2}+\|\chi_\eps(x+x_0)\pd_xu_n\|_{L^2}
\le C(n,\eps)\|\chi_{\eps}(x+x_0)f_n\|_{L^2}\,.
$$
Combining the above with \eqref{eq:Lfn}, we have 
$$
\|\chi_\eps(x+x_0)\pd_x^{-1}u_n\|_{L^2}\le \frac{2C(n,\eps)+1}{n}\|\chi_{\eps}(x+x_0)f_n\|_{L^2}\quad\text{for $n\ne0$.}
$$
To prove \eqref{eq:fred3},
it suffices to show that $\sup_{n\ne0}C(n,\eps)<\infty$.
Let $\tilde{u}_n(x)=\chi_\eps(x+x_0)u_n(x)$,
$\tilde{g}_n(x)=\chi_\eps(x+x_0)(in-\pd_x^2)^{-1}f_n'(x)$ and
$$
\widetilde{T}_n(u)(x)=\chi_\eps(x+x_0)\int_\R \frac{e^{-\sqrt{in}|x-t|}Q(t)u(t)}{\chi_\eps(t+x_0)}
\operatorname{sgn}(t-x)dt,
$$
where $\sqrt{in}$ is chosen with positive real part. The definition of $\widetilde{T}_n$ is obtained thanks to the explicit formula of the kernel
$(in-\pd_x^2)^{-1}\delta=\frac{1}{2\sqrt{in}}e^{-\sqrt{in}|x|}$.
Then \eqref{eq:Lfn} can be rewritten as
$$
\tilde{u}_n=\widetilde{T}_n\tilde{u}_n+\tilde{g}_n\,.$$ 
By \eqref{eq:chch-1},
\begin{equation*}
\Big|\frac{\chi_\eps(x+x_0)}{\chi_\eps(t+x_0)}e^{-\sqrt{in}|x-t|}\Big|\le 2e^{-(\sqrt{|n|/2}-2\eps)|x-t|},
\end{equation*}
and it follows from a convolution estimate that there exists an $n_0\in\N$ such that
$\|\widetilde{T}_n\|_{B(L^2(\R))}\le \frac12$ and
$\|\tilde{g}_n\|_{L^2}\le \frac12\|\chi_\eps(x+x_0)f_n\|_{L^2}$
for $|n|\ge n_0$. Thus we have
$$\|\chi_\eps(x+x_0)u_n\|_{L^2(\R)}\le C\|\chi_\eps(x+x_0)f_n\|_{L^2(\R)}
\quad\text{for $|n|\ge n_0$.}$$
We shall estimate $\chi_{\eps}(x+x_0)\pd_x u_n(x)$ in $L^2$ by a similar argument. We can write
$$
\pd_xu_n(x)=\int_\R e^{-\sqrt{in}|x-t|}(Q(t)u_n(t))'\operatorname{sgn}(t-x)dt+(-\partial_x^2+in)^{-1}f_{n}''(x).
$$
Using a convolution estimate, as above, we obtain that there exist $C>0$ and  $n_0\in\N$ such that for $|n|\geq n_0$,
\begin{multline*}
\Big\|
\chi_\eps(x+x_0)\int_\R e^{-\sqrt{in}|x-t|}(Q(t)u_n(t))'\operatorname{sgn}(t-x)dt\Big\|_{L^2}
\\
\leq C|n|^{-\frac{1}{2}}\big(\|\chi_\eps(x+x_0)u_n\|_{L^2}+\|\chi_\eps(x+x_0)\pd_{x} u_n\|_{L^2}\big).
\end{multline*}
Similarly, we obtain that there exists a constant $C$, independent of $n$ such that for $|n|\geq n_0$,
$$
\|\chi_{\eps}(x+x_0)(-\partial_x^2+in)^{-1}f_{n}''(x)\|_{L^2}\leq C\|\chi_\eps(x+x_0)f_n\|_{L^2}\,.
$$
Therefore, we obtain that there exists a positive constant $C$ and $n_0\in\N$
such that for $|n|\ge n_0$, 
$$\|\chi_\eps(x+x_0)\pd_xu_n\|_{L^2(\R)}\le C\|\chi_\eps(x+x_0)f_n\|_{L^2(\R)}.$$
Thus we have $\sup_{n\ne0}C(n,\eps)<\infty$.
This completes the proof of Lemma \ref{lem:heatop3}.
\end{proof}
\section{Stability of kink solutions of mKP-II}
\label{sec:kink}
In this section, we prove a stability property of the kink solutions 
of the mKP-II equation which will be of crucial importance in the proof
of stability of the line soliton of the KP-II equation.
Recall that $Z=\{u\in H^{8}\,|\, \pd_x^{-1}\pd_yu, \, \pd_xu \in H^8\}.$
We will prove that kinks are stable if the perturbation belongs to
$Z$ and if it is small in $\mathcal{E}$.
\begin{proposition}\label{prop:kinkstability}
For every $\varepsilon>0$, there exists a $\delta>0$ such that if the initial data \eqref{data} of \eqref{eq:MKPII} satisfies
$\|w_0\|_{\mathcal{E}(\R_x\times\T_y)}< \delta$ and $w_0\in Z$,
then there exists a continuous function $\gamma(t)$ such that for every
$t\in\R$, the corresponding solution of \eqref{eq:MKPII} satisfies
$$
\|v(t,x,y)-Q_c(x+\gamma(t))\|_{\mathcal{E}(\R_x\times\T_y)}<\varepsilon.
$$
\end{proposition}
\begin{proof}
As long as $v(t,x,y)$ stays in a small neighborhood of
$\{Q_c(\cdot+x_0)\,|\, x_0\in\R\}$, we can decompose $v(t,x,y)$ as
\begin{equation}
  \label{eq:decomposition}
v(t,x,y)=Q_c(x+\gamma(t))+w(t,x,y), 
\end{equation}
so that $w(t,x,y)$ satisfies the orthogonality condition
\begin{equation}\label{OC}
(w(t,x,y),Q'_c(x+\gamma(t)))
=\frac12(w(t,x,y),\varphi_c(x+\gamma(t)))=0,
\end{equation}
where $(\cdot,\cdot)$ denotes the $L^2(\R_x\times\T_y)$ scalar product.
Note that $\gamma$ depends continuously on $t$. Let
$$
\mathcal{L}_{c,\gamma(t)}\equiv-\pd_x+\pd_x^{-1}\pd_y-2Q_c(x+\gamma(t))\,.
$$
Let us expand the square of the $L^2$ norm of $M^c_{-}(v)$. 
Recalling that $M^c_{-}(Q_c)=0$, we have
\begin{multline*}
\|M^c_{-}(Q_c(x+\gamma(t))+w(t,x,y))\|_{L^2(\R_x\times\T_y)}^2=
\int_{\R_x\times\T_y} (\mathcal{L}_{c,\gamma(t)}w-w^2)^2dxdy
\\=\int_{\R_x\times\T_y}(\mathcal{L}_{c,\gamma(t)}w)^2dxdy
+\int_{\R_x\times\T_y}(w^4-2w^2\mathcal{L}_{c,\gamma(t)}w)dxdy.
\end{multline*}
Thanks to \eqref{OC}, we see that $w(t,x-\gamma(t),y)$ is orthogonal in
$L^2(\R_x\times\T_y)$ to $\varphi_c(x)$. Therefore,
using Lemma~\ref{lem:heatop},
we obtain that there exists a positive constant
$\nu$, independent of $t$ and $w$ such that 
$$
\int_{\R_x\times\T_y}(\mathcal{L}_{c,\gamma(t)}w)^2dxdy
=
\int_{\R_x\times\T_y}(\mathcal{L}_c( w(t,x-\gamma(t),y)  )^2dxdy
\geq
\nu\|w\|^2_{\mathcal{E}(\R_x\times\T_y)}\,.
$$ 
Next we invoke the anisotropic Sobolev embedding of Lemma~\ref{Besov} to have
\begin{equation}\label{eq:l2expand}
\|M^c_{-}(v)(t,x,y)\|_{L^2(\R_x\times\T_y)}^2
\ge \frac{\nu}{2}\|w(t,\cdot)\|_{\mathcal{E}}^2-C\|w(t,\cdot)\|_{\mathcal{E}}^3\,.
\end{equation}
Thanks to the conservation law of Proposition~\ref{prop:wpmkp}, we have
\begin{equation}\label{CL}
\|M^c_{-}(v)(t,x,y)\|_{L^2(\R_x\times\T_y)}^2
=
\|M^c_{-}(Q_c(x)+w_0(x,y))\|_{L^2(\R_x\times\T_y)}^2\,.
\end{equation}
Now expanding the square of the $L^2$ norm of $M^c_{-}(Q_c(x)+w_0(x,y))$ and
using Lemma~\ref{Besov}, we have
\begin{equation}\label{exbis}
\|M^c_{-}(Q_c(x)+w_0(x,y))\|_{L^2(\R_x\times\T_y)}^2
\leq C(\|w_0\|^2_{\mathcal{E}}+\|w_0\|^4_{\mathcal{E}}).
\end{equation}
Combining \eqref{eq:l2expand}, \eqref{CL} and \eqref{exbis}, we get
$\|w(t,\cdot)\|_{\mathcal{E}}\leq C\|w_0\|_{\mathcal{E}}$ 
provided $\delta\ll 1$.
This completes the proof of Proposition~\ref{prop:kinkstability}.
\end{proof}
\section{Asymptotic stability of kink solutions of mKP-II}
\label{sec:asympst}
In this section, we will prove asymptotic stability of kink solutions.
\begin{proposition}
  \label{prop:asympkink}
There exists a $\delta>0$ such that if the initial data \eqref{data} of
\eqref{eq:MKPII} satisfies $\|w_0\|_{\mathcal{E}(\R_x\times\T_y)}< \delta$ and $w_0\in Z$,
then the corresponding solution of \eqref{eq:MKPII} satisfies
$$
\lim_{t\to\infty}\|v(t,x,y)-Q_c(x+\gamma(t))\|_{\mathcal{E}((x\ge-2ct) \times\T_y)}=0,
$$
where $\gamma(t)$ is a $C^1$-function satisfying
$\dot{\gamma}(t)=c+O(\|w_0\|_{\mathcal{E}})$.
\end{proposition}
To prove Proposition \ref{prop:asympkink}, we will use the Miura transform
$u(t,x,y)=M_-^c(v)(t,x-3ct,y)$ and the monotonicity property of 
the KP-II equation that follows from the Kato type estimate for the KP-II
equation. Let us recall the Kato type identity for the KP-II equation.
\begin{lemma}{\cite[Lemma 1]{dBM}}
  \label{lem:kato}
Let $u(t)\in C(\R;H^8(\R_x\times\T_y))$ be a solution of \eqref{eq:KPII} and
$\phi(x)$ be of the class $C^3$. Then
\begin{equation}
  \label{eq:kato}
  \begin{split}
\frac{d}{dt}\int_{\R_x\times \T_y} u^2(t,x,y)\phi(x)dxdy= &
\int_{\R_x\times \T_y}\left(-3(\pd_xu)^2-3(\pd_x^{-1}\pd_yu)^2-4u^3\right)
\phi'dxdy
\\ & +\int_{\R_x\times \T_y}u^2\phi'''dxdy.
  \end{split}
\end{equation}
\end{lemma}
Next, we will show that small solutions of
the KP-II equation locally tends to $0$ as $t\to\infty$ by using
Lemma \ref{lem:kato}.
\begin{lemma}\label{lem:monotonicity}
Let again $\chi_\eps(x)=(1+\tanh (\eps x))/2$ and let $u(t)\in C(\R;H^8(\R_x\times\T_y))$ be
a solution of the KP-II equation  \eqref{eq:KPII}.
For every $c_1>0$, there exist positive numbers $\delta$ and $\eps_0$ such
that if $\eps\in(0,\eps_0)$ and $\|u(0,x,y)\|_{L^2(\R_x\times \T_y)}\le \delta$,
then for every $x_0\in\R$,
\begin{gather*}
\frac{d}{dt}\int_{\R_x\times \T_y} u^2(t,x,y)\chi_\eps(x+x_0-c_1t)dxdy\le 0,\\
\lim_{t\to\infty}\int_{\R_x\times \T_y} u^2(t,x,y)\chi_\eps(x-c_1t)dxdy=0.
\end{gather*}
\end{lemma}
To show Lemma \ref{lem:monotonicity}, we use the following.
\begin{claim}\label{amiens}
Let $\eps>0$. 
There exists a positive constant $C$ such that for every $w\in\mathcal{E}(\R_x\times\T_y)$ and $x_0\in\R$, 
\begin{equation}
\label{eq:ias}
\Big(\int_{\R_x\times\T_y}(\chi'_\eps(x+x_0))^2w^4dxdy\Big)^{\frac 1 2}\le 
C\int_{\R_x\times \T_y}\Big((\pd_x w)^2+(\pd_x^{-1}\pd_y w)^2+w^2\Big)
\chi_\eps'(x+x_0)dxdy.
\end{equation}
\end{claim}
\begin{proof}
This claim can be proved in a similar way as  \cite[Lemma 2]{MST07}.
By a density argument, we can suppose that $w\in Z$ and in particular $w$
vanishes as $x\rightarrow\pm\infty$.
\par 
Suppose that $w$ vanishes at a point of $\T_y$, say $y=0$.
Using Fubini's theorem and integration by parts, we have 
\begin{equation}
  \label{eq:*1}
  \begin{split}
\int_\R \chi'_\eps(x+x_0)w^2(x,y)dx
=& 2\int_\R\chi'_\eps(x+x_0)\Big(\int_{0}^yw(x,t)w_y(x,t)dt\Big)dx
\\=& -2\int_{0}^y
\Big(\int_\R\pd_x(\chi'_\eps(x+x_0)w(x,t))\,  \pd_x^{-1}\pd_yw(x,t)dx\Big)dt
\\ \le & C
\int_{\R_x\times \T_y}\Big((\pd_x w)^2+(\pd_x^{-1}\pd_y w)^2+w^2\Big)
\chi_\eps'(x+x_0)dxdy,    
  \end{split}
\end{equation}
where we used an inequality $|\chi''_\eps|\leq 4\eps \chi'_\eps$
and the inequality $2ab\leq a^2+b^2$.
If $w(x,y)$ does not vanish on $\T_y$, we apply a partition of unity argument
to $w$. More precisely, let $\psi_1(y)$ and $\psi_2(y)$ be nonnegative
smooth functions such that $\psi_1+\psi_2\equiv1$ and $\psi_1(y)=0$
for $|y|\le \pi/5$ and $\psi_2(y)=0$ for $y\notin (-\pi/4,\pi/4)$. Then
\begin{align*}
& \int_\R \chi'_\eps(x+x_0)\psi_1(y)w^2(x,y)dx
\\=& -2\int_{0}^y\psi_1(t)
\Big(\int_\R\pd_x(\chi'_\eps(x+x_0)w(x,t))\pd_x^{-1}\pd_yw(x,t)dx\Big)dt
\\& +\int_{0}^y\psi_1(t)\Big(\int_\R\chi'_\eps(x+x_0)w^2(x,t)dx\Big)dt
\\ \leq C & 
\int_{\R_x\times \T_y}
\Big((\pd_x w)^2+(\pd_x^{-1}\pd_y w)^2+w^2\Big)\chi_\eps'(x+x_0)dxdy.
\end{align*}
We can estimate $\int_\R \chi'_\eps(x+x_0)\psi_2(y)w^2(x,y)dx$ in the same way, by replacing $0$ with another point on $\T_y$, say $\pi/3$.
Thus we have 
\begin{equation}\label{isa2}
\sup_{y}\int_\R \chi'_\eps(x+x_0)w^2(x,y)dx \le  C
\int_{\R_x\times \T_y}\Big((\pd_x w)^2+(\pd_x^{-1}\pd_y w)^2+w^2\Big)
\chi_\eps'(x+x_0)dxdy.
\end{equation}
Write
$$\chi'_{\eps}(x+x_0)w^2(x,y)=\int_{-\infty}^{x}
\big(\chi''_{\eps}(z+x_0)w^2(z,y)+2\chi'_{\eps}(z+x_0)
\partial_{x}w(z,y)w(z,y)\big)dz\,.$$
This yields
\begin{equation}
  \label{eq:ias3}
\sup_x\big[\chi'_{\eps}(x+x_0)w^2(x,y)\big] \leq C
\int_{-\infty}^{\infty}(w^2(z,y)+(\partial_x w)^2(z,y))\chi'_{\eps}(z+x_0)dz\,.  
\end{equation}
Applying \eqref{eq:ias3} to the first factor of
$(\chi'_\eps(x+x_0))^2w^4(x,y)=(\chi'_\eps(x+x_0)w^2)\times(\chi'_\eps(x+x_0)w^2)$
and integrating the resulting inequality over $\R_x\times\T_y$, we obtain
\begin{multline}\label{isa1}
\int_{\R_x\times\T_y}(\chi'_\eps(x+x_0))^2w^4(x,y)dxdy
\leq 
\\
C\,\int_{-\infty}^{\infty}
\Big(\int_{-\infty}^{\infty}w^2(x,y)\chi'_{\eps}(x+x_0)dx\Big)
\Big(\int_{-\infty}^{\infty}
(w^2(z,y)+(\partial_x w)^2(z,y))\chi'_{\eps}(z+x_0)dz\Big)
dy.
\end{multline}
Substituting \eqref{isa2} into \eqref{isa1}, we obtain \eqref{eq:ias}.
This completes the proof of Claim~\ref{amiens}.
\end{proof}
\begin{proof}[Proof of Lemma~\ref{lem:monotonicity}]
By taking into account the time translation and using Lemma~\ref{lem:kato}, we car write 
\begin{equation}
  \label{eq:kato2}
  \begin{split}
& \frac{d}{dt}\int_{\R_x\times \T_y} u^2(t,x,y)\chi_\eps(x+x_0-c_1t)dxdy
\\=& \int_{\R_x\times \T_y}\left(-3(\pd_xu)^2-3(\pd_x^{-1}\pd_yu)^2-4u^3\right)
(t,x,y)\chi_\eps'(x+x_0-c_1t)dxdy
\\ & +\int_{\R_x\times \T_y}u^2(t,x,y)(\chi_\eps'''(x+x_0-c_1t)-c_1\chi_\eps'(x+x_0-c_1t))dxdy.    
  \end{split}
\end{equation}
Using Claim~\ref{amiens} , we have
\begin{multline}\label{isa3}
\Big|\int u^3(t,x,y)\chi_\eps'(x+x_0-c_1t)dxdy\Big|
\leq \|u(t,\cdot)\|_{L^2}
\Big(\int(\chi'_\eps(x+x_0-c_1 t))^2u^4dxdy\Big)^{\frac 1 2}
\\
\leq 
C\|u(0,\cdot)\|_{L^2}
\int_{\R_x\times \T_y}\Big((\pd_x u)^2+(\pd_x^{-1}\pd_y u)^2+u^2\Big)(t,x,y)\chi_\eps'(x+x_0-c_1 t)dxdy\,.
\end{multline}
It follows from \eqref{eq:kato2}, \eqref{isa3} and the fact that  $|\chi_\eps'''|\le 8\eps^2\chi_\eps'$ that
\begin{align*}
& \frac{d}{dt}\int_{\R_x\times \T_y} u^2(t,x,y)\chi_\eps(x-c_1t)dxdy
\\ \leq &
-C\,\int_{\R_x\times \T_y}\left((\pd_xu)^2+(\pd_x^{-1}\pd_yu)^2+u^2\right)(t,x,y)
\chi_\eps'(x-c_1t)dxdy    
\end{align*}
provided $\epsilon$ and $\delta$ are sufficiently small.
Thus we have the former part of Lemma \ref{lem:monotonicity}.
\par
Let $$I_{x_0}(t)=\int_{\R_x\times\T_y}u^2(t,x,y)\chi_\eps(x-\frac12c_1t-x_0)dxdy.$$
Taking $\delta>0$ smaller if necessary, we have
$I_{x_0}(t)\le I_{x_0}(0)$ for every $t\ge 0$ and $x_0\geq 0$.
Since $\lim_{x_0\to\infty} I_{x_0}(0)=0$, we have
$$
\lim_{t\to\infty}\int_{\R_x\times\T_y}u^2(t,x,y)\chi_\eps(x-c_1t)=\lim_{t\to\infty}I_{c_1t/2}(t)\leq \liminf_{t\to\infty}I_{c_1t/2}(0)=0.
$$
Thus we complete the proof of Lemma~\ref{lem:monotonicity}.
\end{proof}
Next, we will derive a modulation equation on the parameter
$\gamma(t)$ which describes the phase shift of the line soliton.
\begin{lemma}\label{lem:modeq}
Let $w_0\in Z$ and $v(t,x,y)$ be a solution of \eqref{eq:MKPII}
satisfying \eqref{data}. Let $\gamma(t)$ be a $C^1$-function 
satisfying \eqref{eq:decomposition} and \eqref{OC}.
There exist positive constants $\delta>0$ and $C$ such that
if $\|w_0\|_{\mathcal{E}}<\delta$,
$|\dot{\gamma}(t)-c|+|\gamma(0)|\le C\|w_0\|_{\mathcal{E}}$.
\end{lemma}
\begin{proof}
 By Proposition \ref{prop:kinkstability}, a solution $v(t)$ remains in a 
tubular neighborhood of the kink solutions for every $t\in\R$ if
$\delta$ is sufficiently small. Applying the implicit function theorem
and using a continuation argument, we find a $C^1$-function $\gamma(t)$ that
satisfies $|\gamma(0)|=O(\|w_0\|_{\mathcal{E}})$ and
\eqref{eq:decomposition} and \eqref{OC} for every $t\in\R$.
\par
Differentiating \eqref{OC} with respect to $t$ and substituting 
\eqref{eq:MKPII} into the resulting equation, we have
\begin{align*}
& \frac{d}{dt}\int w(t,x,y)Q_c'(x+\gamma(t))dxdy
\\=&(c-\dot{\gamma}(t))\{\|Q_c'\|_{L^2}^2-(w(t),Q_c''(x+\gamma(t)))\}
+O(\|w(t)\|_{\mathcal{E}})
\\=&0.
\end{align*}
Thus we complete the proof of Lemma~\ref{lem:modeq}.
\end{proof}
We will make use of the following statement 
(which will be used in the particular case $p=4$).
\begin{claim}\label{cl:weightimbed}
Let $2\le p\le 6$ and $\eps>0$. There exists a positive constant $C$ such that
for every $w\in\mathcal{E}(\R_x\times\T_y)$ and $x_0\in\R$,
\begin{equation}\label{eq:wimbed}
\int_{\R_x\times\T_y}\chi_\eps(x+x_0)|w(x,y)|^pdxdy
\le C\|w\|_{\mathcal{E}(\R_x\times\T_y)}^{p-2}\|w\|_{\mathcal{E}_{\eps,x_0}}^2.
\end{equation}
\end{claim}
\begin{proof}
The proof follows the line of the proof of Claim~\ref{amiens}.
We only need to consider the endpoint statements, then the other 
cases follow by the H\"older inequality with respect to the measure
$\chi_{\eps}(x+x_0)dxdy$.
The case $p=2$ is obvious. Let us consider the case $p=6$.
First, by invoking the inequality $|\chi'_{\eps}|\leq 2\eps \chi_\eps$
in place of $|\chi''_\eps|\leq 4\eps \chi'_\eps$ in the proof of \eqref{isa2},
we get that there exists a constant $C$ such that for every $x_0\in\R$, every
$y\in\T$ and every $w\in \mathcal{E}$,
\begin{equation}\label{isa4}
\sup_{y}\int_{-\infty}^{\infty} \chi_\eps(x+x_0)|w(x,y)|^2dx\leq
C\|w\|_{\mathcal{E}_{\eps,x_0}}^2\,.
\end{equation}
Similarly, we have
\begin{equation}\label{isa5}
\sup_{y}\int_{-\infty}^{\infty} |w(x,y)|^2dx\leq   C\|w\|_{\mathcal{E}}^2.
\end{equation}
Combining \eqref{isa4} and \eqref{isa5} with
$\sup_x|w(x,y)|^4\le 2(\int w(z,y)^2dz)(\int w_x(z,y)^2dz)$, we obtain
\begin{align*}
& \int_{\R_x\times\T_y} \chi_\eps(x+x_0)|w(x,y)|^6dxdy \\ \le &
2\int_{-\infty}^{\infty}\left\{\Big(\int_\R |w(z,y)|^2dz\Big)\Big(\int_\R |w_x(z,y)|^2dz\Big)\int_\R\chi_\eps(x+x_0)|w(x,y)|^2dx\right\}dy
\\ \le & C\|w\|_{\mathcal{E}(\R_x\times\T_y)}^4\|w\|_{\mathcal{E}_{\eps,x_0}}^2.
\end{align*}
Thus we complete the proof of Claim~\ref{cl:weightimbed} .
\end{proof}
Now we are in position to prove Proposition \ref{prop:asympkink}.
\begin{proof}[Proof of Proposition \ref{prop:asympkink}]
Let $u_-(t,x,y)=M_-^c(v)(t,x-3ct,y)$. Then  $u_-(t,x,y)\in C(\R;H^8(\R_x\times\T_y))$
and $u_-(t,x,y)$ is a solution of \eqref{eq:KPII}.
Lemma \ref{lem:monotonicity} implies that there exist positive numbers
$\eps$ and $\delta$ such that if $\|w_0\|_{\mathcal{E}}<\delta$,
\begin{equation}
  \label{eq:locconv1}
\lim_{t\to\infty}\int \chi_\eps(x+2ct)|M_-^c(v)(t,x,y)|^2dxdy=0.  
\end{equation}
Let us decompose $v(t,x,y)$ as \eqref{eq:decomposition} and \eqref{OC}. Thanks to the imposed orthogonality conditions,
Lemma \ref{lem:heatop3} implies 
\begin{equation}\label{eq:lwb2}
\|\mathcal\mathstrut{\sqrt{\chi_\eps}}(x+2ct)
(\mathcal{L}_{c,\gamma(t)}w)(t,x,y)\|_{L^2}\ge \nu \|w(t,x,y)\|_{\mathcal{E}_{\eps,2ct}},
\end{equation}
where $\nu$ is a positive constant depending only on $c$. Using \eqref{eq:lwb2}, the Cauchy-Schwarz inequality and Claim~\ref{cl:weightimbed} with $p=4$,
we have
\begin{multline}\label{eq:lwb1}
\|\mathstrut{\sqrt{\chi_\eps}}(x+2ct)M_-^c(v)(t,x,y)\|_{L^2} \ge  
\frac{\nu}{2} \|w(t,x,y)\|_{\mathcal{E}_{\eps,2ct}}
\\
-C\Big(\int_{\R_x\times\T_y}w^4(t,x,y)\chi_\eps(x+2ct)\Big)^{\frac{1}{2}}
\geq 
\frac{\nu}{2} \|w(t,x,y)\|_{\mathcal{E}_{\eps,2ct}}-C\|w\|_{\mathcal{E}}\|w\|_{\mathcal{E}_{\eps,2ct}}.
\end{multline}
Since $\sup_{t\in\R}\|w(t,\cdot)\|_{\mathcal{E}}=O(\delta)$ by
Proposition~\ref{prop:kinkstability},  coming back to \eqref{eq:locconv1} ,
we get
$$
\lim_{t\to\infty}\|w(t,\cdot)\|_{\mathcal{E}_{\eps,2ct}}=0.
$$
This completes the proof of Proposition \ref{prop:asympkink}.
\end{proof}
\section{Proof of Theorem \ref{thm1}}
\label{sec:proof}
We are now in position to prove Theorem~\ref{thm1}.
\begin{proof}[Proof of Theorem \ref{thm1}]
Let $u_0(x,y)=\varphi_c(x)+\tilde{u}_0(x,y)$ and let $u(t,x,y)$ be the solution of \eqref{eq:KPII} with data $u_0$.
Fix $\varepsilon>0$.
By Proposition~\ref{prop:bijMiura}, for every $\varepsilon_1>0$ there exists $\delta_1>0$ such that if $\|\tilde{u}_0\|_{L^2}< \delta_1$,
there exists a unique couple $(k,w_0)\in\R\times\mathcal{E}$ such that
\begin{equation}\label{IF}
|k-c|<\varepsilon_1, \quad\|w_0\|_{\mathcal{E}}<\varepsilon_1,\quad
M^k_{+}(Q_k+w_0)=\varphi_c+\tilde{u}_0.
\end{equation}
Let 
$w_{0,n}\in Y$ $(n\in\N)$ be such that 
\begin{equation}\label{conv}
\lim_{n\to\infty}\|w_{0,n}-w_0\|_{\mathcal{E}}=0
\end{equation}
and let
$v_n(t)$ be a solution of \eqref{eq:MKPII} with initial data
$v_n(0)=Q_k+w_{0,n}$. 
Proposition~\ref{prop:kinkstability}, \eqref{IF} and \eqref{conv} imply that
for any $\varepsilon_2>0$, there exist $\delta_2>0$ and $n_0\in\N$ such that 
if $\|\tilde{u}_0\|_{L^2}<\delta_2$ and $n\geq n_0$,
then there exist continuous functions $\tilde{\gamma}_n(t)$ such that
\begin{equation}\label{eq:vnstability}
\sup_{t\in\R}\,\sup_{n\geq n_0}
\|v_n(t,x,y)-Q_k(x+\tilde{\gamma}_n(t))\|_{\mathcal{E}(\R_x\times\T_y)}<\varepsilon_2.
\end{equation}
Let
$u_n(t,x,y)\equiv M^k_{+}(v_n)(t,x-3kt,y).$
Then $u_n(t,x,y)$ is a solution of the KP-II equation \eqref{eq:KPII}.
Thanks to \eqref{conv}, the sequence $u_{n}(0)=M_k^{+}(Q_k+w_{0,n})$ converges, as $n\rightarrow\infty$,
in $L^2(\R_x\times\T_y)$
to $M_k^{+}(Q_k+w_0)=\varphi_c+\tilde{u}_0$. Therefore
\begin{equation}\label{conv_data}
\lim_{n\to\infty}\|u_n(0,x,y)-\varphi_c(x)-\tilde{u}_0(x,y)\|_{L^2(\R_x\times\T_y)}=0.
\end{equation}
(recall that $u(0,x,y)=\varphi_c(x)+\tilde{u}_0(x,y)$).
Next, using Lemma~\ref{Besov}, we have
\begin{multline*}
\|u_n(t,x,y)-\varphi_k(x+\gamma_n(t))\|_{L^2(\R_x\times\T_y)}
\\
=
\|M^k_{+}(v_n)-M^k_{+}(Q_{k}(x+\tilde{\gamma}_n(t)))\|_{L^2(\R_x\times\T_y)}
\\
\leq
C\|v_n(t,x,y)-Q_k(x+\tilde{\gamma}_n(t))\|_{\mathcal{E}(\R_x\times\T_y)}
\\
\times
\big(1+\|v_n(t,x,y)\|_{\mathcal{E}(\R_x\times\T_y)}+\|Q_k\|_{L^\infty(\R_x\times\T_y)}\big)\,,
\end{multline*}
where $\gamma_n(t)=\tilde{\gamma}_n(t)-3kt$.
Coming back to \eqref{eq:vnstability}, we obtain that
for any $\varepsilon_3>0$, there exist $\delta_3>0$ and $n_0\in\N$ such that 
if $\|\tilde{u}_0\|_{L^2}<\delta_3$,
\begin{equation}\label{pak}
\|u_n(t,x,y)-\varphi_k(x+\gamma_n(t))\|_{L^2(\R_x\times\T_y)}<\varepsilon_3,\quad
\forall\,t\in\R,\quad \forall\, n\geq n_0.
\end{equation}
Using the triangle inequality, we obtain
\begin{multline*}
\inf_{\gamma\in\R}\|u(t,x,y)-\varphi_c(x+\gamma)\|_{L^2(\R_x\times\T_y)}
\\
\leq
\|u(t,x,y)-\varphi_c(x+\gamma_n(t))\|_{L^2(\R_x\times\T_y)}
\leq
\|u(t,\cdot)-u_n(t,\cdot)\|_{L^2(\R_x\times\T_y)}
\\
+
\|u_n(t,x,y)-\varphi_k(x+\gamma_n(t))\|_{L^2(\R_x\times\T_y)}
+
\|\varphi_k-\varphi_c\|_{L^2(\R_x\times\T_y)}\,.
\end{multline*}
For any $t\in\R$, using \eqref{conv_data} and the $L^2$ well-posedness
result of \cite{MST}, we see that there exists $n_1$ (depending on $t$ and $\varepsilon$) such that
for $n\geq n_1$,
$$\|u(t,\cdot)-u_n(t,\cdot)\|_{L^2(\R_x\times\T_y)}<\frac{\varepsilon}{3}\,.
$$
By \eqref{pak},  there exists $n_2$ such that for $\|\tilde{u}_0\|_{L^2}\ll 1$
and $n\geq n_2$,
$$
\|u_n(t,x,y)-\varphi_k(x+\gamma_n(t))\|_{L^2(\R_x\times\T_y)}<\frac{\varepsilon}{3}.
$$
Finally, using \eqref{IF}, we obtain that for $\|\tilde{u}_0\|_{L^2}\ll 1$,
$$
\|\varphi_k-\varphi_c\|_{L^2(\R_x\times\T_y)}<\frac{\varepsilon}{3}.
$$
Summarizing, we obtain that for $\|\tilde{u}_0\|_{L^2}\ll 1$,
$$
\inf_{\gamma\in\R}\|u(t,x,y)-\varphi_c(x+\gamma)\|_{L^2(\R_x\times\T_y)}<\varepsilon,\quad \forall\, t\in\R.
$$
This completes the proof of the orbital stability.
\par
Next, we will prove asymptotic stability of line solitons.
Let $a$ and $\delta_4$ be small positive numbers and $n_0$ be a large integer.
Then if $\|w_0\|_{\mathcal{E}}\le \delta_4$ and $n\ge n_0$,
\begin{multline*}
\|u_n(t,x,y)-\varphi_k(x+\gamma_n(t))\|_{L^2((x>ct)\times \T_y)}
\\=
\|M_+^k(v_n)(t,x,y)-M_+^k(Q_k(x+\tilde{\gamma}_n(t)))\|_{L^2((x>(c-3k)t)\times \T_y)}
\\ \le  C_1\|v_n(t,x,y)-Q_k(x+\tilde{\gamma}_n(t))\|_{\mathcal{E}_{a,(3k-c)t}},
\end{multline*}
and we obtain 
$$
\|v_n(t,x,y)-Q_k(x+\tilde{\gamma}_n(t))\|_{\mathcal{E}_{a,(3k-c)t}}\le 
C_2
\|M_-^k(v_n)(t,x,y)\|_{L^2_{a,(3k-c)t}}
$$
in the same way as \eqref{eq:lwb1}.
We remark that $C_1$ and $C_2$ do not depend on $n$.
On the other hand, we can write
$$
\|M_-^k(v_n)(t,x,y)\|_{L^2_{a,(3k-c)t}}^2=
\int_{\R_x\times \T_y}\chi_\eps(x-\frac{ct}{2}-\frac{ct}{2})\big(M_-^k(v_n)(t,x-3kt,y)\big)^2dxdy
$$
and using Lemma~\ref{lem:monotonicity} with $c_1=c/2$ and $x_0=-ct/2$,
we obtain that
\begin{align*}
& \|M_-^k(v_n)(t,x,y)\|_{L^2_{a,(3k-c)t}}^2
\\ \le & \int_{\R_x\times \T_y}\chi_\eps(x-\frac{ct}{2})\big(M_-^k(Q_k(x)+w_{0,n}(x,y))\big)^2dxdy
\\ \leq & 
C\int_{\R_x\times \T_y}\chi_\eps(x-\frac{ct}{2})\left\{(\mathcal{L}_k w_0)^2(x,y)
+w_0^4(x,y)\right\}dxdy 
+C\|w_0-w_{0,n}\|_{\mathcal{E}}^2,
\end{align*}
where $C$ is a constant that does not depend on $n$.
Therefore using the dominated convergence,
we obtain that there exist sequences $\{\eps_j\}_{j\ge1}$ and $\{n_0(j)\}_{j\ge1}$
such that $\lim_{j\to\infty}\eps_j=0$ and 
$$
\sup_{n\ge n_0(j)}\sup_{t\ge j}
\|u_n(t,x,y)-\varphi_k(x+\gamma_n(t))\|_{L^2((x>ct)\times \T_y)}
<\eps_j.
$$
\par
On the other hand, there exists $n_1(j)\ge n_0(j)$ such that
$$
\sup_{n\ge n_1(j)}\sup_{t\in[j,j+1]}\|u(t)-u_n(t)\|_{L^2(\R_x\times\T_y)}<
\eps_j
$$
thanks to the well-posedness of \eqref{eq:KPII} for data in $L^2(\R_x\times\T_y)$
(see \cite{MST}). Letting $\tilde{c}=k$ and 
$x(t)=-\gamma_{n_1(j)}(t)$ for $t\in[j,j+1]$, we have \eqref{eq:asympst}.
In view of  Lemma \ref{lem:modeq}, we have
$$|\dot{\gamma}_n(t)+2c|+|\gamma_n(0)| \le C\|w_{0,n}\|_{\mathcal{E}},$$
where $C$ is a constant that does not depend on $n$ and $t$.
This completes the proof of Theorem~\ref{thm1}.
\end{proof}
\section*{Acknowledgment}
The main part of this work was carried out when the first author stayed
at Universit\'e de Cergy-Pontoise.
He would like to express his gratitude to for its hospitality.


\begin{thebibliography}{10}
\bibitem{Benjamin} T.~Benjamin,
 {\it The stability of solitary waves},
 Proc. R. Soc. Lond. A 328 (1972), 153-183.
%
\bibitem{BIN} O.~Besov, V.~Ilin and S.~Nikolski,
{\it Integral representations of functions and embedding theorems},
J. Wiley 1978. 
%
\bibitem{Bona} J.~L.~Bona,
{\it The stability of solitary waves},
Proc. R. Soc. Lond. A 344 (1975), 363-374.
%
\bibitem{dBM} A.~de Bouard and Y.~Martel,
{\it Non existence of $L^2$-compact solutions of the Kadomtsev-Petviashvili
II equation}, Math. Ann. 328 (2004) 525-544.
%
\bibitem {Bourgain} J.~Bourgain,
{\it On the Cauchy problem for the Kadomtsev-Petviashvili equation},
 GAFA 3 (1993), 315-341.
%
\bibitem{Cu} S.~Cuccagna,
{\it On asymptotic stability in 3D of kinks for the $\phi^4$ model},
Trans. Amer. Math. Soc. (2008), 2581-2614. 
%
\bibitem{FP} G.~Friesecke and R.~L.~Pego, 
{\it Solitary waves on FPU lattices. I, Qualitative properties, renormalization and continuum limit}, Nonlinearity 12 (1999), 1601--1627;
{\it Solitary waves on FPU lattices. II, Linear implies nonlinear stability},
Nonlinearity 15 (2002), 1343--1359;
{\it Solitary waves on Fermi-Pasta-Ulam lattices.
III, Howland-type Floquet theory}, Nonlinearity 17 (2004), 207--227;
{\it
Solitary waves on Fermi-Pasta-Ulam lattices. IV,
Proof of stability at low energy},
Nonlinearity 17 (2004), 229--251.
%
\bibitem{GPS} A.~Gr\"unrock, M.~Panthee and J.~Drumond Silva,
{\it On KP-II equations on cylinders}, Ann. IHP Analyse non lin\'eaire 26 (2009), 2335-2358.
%
\bibitem{Hadac} M.~Hadac,
{\it Well-posedness of the KP-II equation and generalizations}, Trans. Amer. Math. Soc. (2008), 6555-6572.
%
\bibitem{HHK} 
M.~Hadac and S.~Herr and H.~Koch, 
{\it Well-posedness and scattering for the KP-II equation in a critical space},
Ann. IHP Analyse non lin\'eaire 26 (2009), 917-941.
%
\bibitem{IKT}
A.~Ionescu, C.~Kenig and D.~Tataru,
{\it Global well-posedness of the initial value problem for the initial value problem for the the KP-I equation in the energy space},
Invent. Math. 173 (2008) 265-304.
%
\bibitem{IM} P.~Isaza and J.~Mejia,
{\it Local and global Cauchy problems for the Kadomtsev-Petviashvili (KP-II)
equation in Sobolev spaces of negative indices},
Comm. Partial Differential Equations 26 (2001), 1027-1057.
%
\bibitem{KMa} C.~Kenig and Y.~Martel,
{\it Global well-posedness in the energy space for a modified KP II equation
via the Miura transform}, 
Trans. Amer. Math. Soc. 358 (2006), no. 6, 2447-2488.
%
\bibitem{KP} B.~B.~Kadomtsev and V.~I.~Petviashvili, 
{\it On the stability of solitary waves in weakly dispersive media}, 
Sov. Phys. Dokl. 15 (1970), 539-541.
%
\bibitem{MM} Y.~Martel and F.~Merle,
{\it A Liouville theorem for the critical generalized Korteweg-de Vries
equation}, J. Math. Pures Appl. 79 (2000),  339-425.
%
\bibitem{MV} F.~Merle, L.~Vega,
{\it $L^2$ stability of solitons for KdV equation},  
Int. Math. Res. Not. (2003), 735-753. 
%
\bibitem{Mi1} T.~Mizumachi,
{\it Asymptotic stability of lattice solitons},
Comm. Math. Phys. 288 (2009), 125-144.
%
\bibitem{MP} T.~Mizumachi and R.~L.~Pego,
{\it Asymptotic stability of Toda lattice solitons}, Nonlinearity  21 (2008),
2099--2111.
%
\bibitem{Molinet} L.~Molinet,
{\it On the asymptotic behavior of solutions to the (generalized) Kadomtsev-Petviashvili-Burgers equation},
J. Diff. Eq. 152 (1999) 30-74.
%
\bibitem{MST07} L.~Molinet, J.~C.~Saut and N.~Tzvetkov,
{\it Global well-posedness for the KP-I equation on the background of
a non localized solution},
Comm. Math. Phys, 272 (2007), 775-810.
%
\bibitem{MST_contr} L. Molinet, J.-C. Saut and N. Tzvetkov, 
{\it Remarks on the mass constraint for KP-type equations},
SIAM J. Math. Anal. 39 (2007), 627-641.
%
\bibitem{MST} L.~Molinet, J.~C.~Saut and N.~Tzvetkov, 
{\it Global well-posedness for the KP-II equation on the background of a non localized solution},
Preprint~2010.
%
\bibitem{RT1} F.~Rousset and N.~Tzvetkov,
{\it Transverse nonlinear instability for two-dimensional dispersive models },
Ann. IHP, Analyse Non Lin\'eaire 26 (2009), 477-496.
%
\bibitem{RT2} F.~Rousset and N.~Tzvetkov, {\it Transverse nonlinear instability for some Hamiltonian PDE's},
J. Math. Pures Appl. 90 (2008), 550-590.
%
\bibitem {Tak} H.~Takaoka, {\it Global well-posedness for the Kadomtsev-Petviashvili II equation}, Discrete Contin. Dynam. Systems, 6 (2000), 483-499.
%
\bibitem{TT} H.~Takaoka and N.~Tzvetkov,
{\it On the local regularity of Kadomtsev-Petviashvili-II equation},
IMRN 8 (2001), 77-114.
%
\bibitem{Tom} M.~Tom, {\it On a generalized Kadomtsev-Petviashvili equation}, Contemporary Mathematics, AMS 200 (1996), 193-210.
%
\bibitem{Tz} N.~Tzvetkov, {\it Global low regularity solutions for Kadomtsev-Petviashvili equation}, Diff. Int. Eq. 13 (2000), 1289-1320.
%
\bibitem{VA} J.~Villarroel and M.~Ablowitz, 
{\it On the initial value problem for the KPII equation with data that do not decay along a line},  
Nonlinearity  17  (2004), 1843-1866.
%
\bibitem{W} M.~V.~Wickerhauser,
{\it Inverse scattering for the heat equation and evolutions in $(2+1)$ variables},
Comm. Math. Phys. 108 (1987), 67-89.
%
\bibitem{Z} V.~Zakharov,
{\it Instability and nonlinear oscillations of solitons},
JEPT Lett. 22 (1975), 172-173.
%
\bibitem{ZS} V.~Zakharov and E.~Schulman, 
{\it Degenerative dispersion laws, motion invariants and kinetic equations},
 Physica D 1 (1980), 192-202.
\end{thebibliography}
\end{document}